\documentclass[11pt,reqno]{amsart}

\usepackage[T1]{fontenc}
\usepackage{lmodern}
\usepackage{microtype}
\usepackage{mathtools}
\usepackage{amssymb}
\usepackage{mathrsfs}
\usepackage{booktabs}
\usepackage{enumitem}
\usepackage[margin=1.12in]{geometry}
\usepackage[pdfusetitle,hidelinks]{hyperref}
\usepackage{aliascnt}
\usepackage[nameinlink,noabbrev,capitalise]{cleveref}
\usepackage{todonotes}

\setlist[enumerate,1]{label=\textup{(\roman*)},leftmargin=2.2em}
\setlist[itemize]{leftmargin=2em}
\numberwithin{equation}{section}

\newtheorem{theorem}{Theorem}[section]
\newaliascnt{proposition}{theorem}
\newtheorem{proposition}[proposition]{Proposition}
\aliascntresetthe{proposition}
\newaliascnt{lemma}{theorem}
\newtheorem{lemma}[lemma]{Lemma}
\newtheorem{question}[theorem]{Question}
\newtheorem{example}[theorem]{Example}

\aliascntresetthe{lemma}
\newaliascnt{corollary}{theorem}
\newtheorem{corollary}[corollary]{Corollary}
\aliascntresetthe{corollary}
\theoremstyle{remark}
\newaliascnt{remark}{theorem}

\aliascntresetthe{remark}

\newcommand{\sshf}[1]{\mathcal{O}_{#1}}
\newcommand{\shf}[1]{\mathscr{#1}}

\newcommand{\prj}[1]{\mathbb{P}^{#1}}

\newcommand{\iso}{\simeq}
\newcommand{\ses}[3]{0\rightarrow#1\rightarrow#2\rightarrow#3\rightarrow{0}}

\newcommand{\paren}[1]{\left(#1\right)}
\newcommand{\is}[1]{\mathscr I_{#1}}
\newcommand{\Sym}{\operatorname{Sym}}
\newcommand{\Pic}{\operatorname{Pic}}

\newcommand{\rk}[1]{\operatorname{rank}(#1)}

\newcommand{\pa}{p_a}
\newcommand{\NP}{\mathrm N}
\newcommand{\floor}[1]{\left\lfloor #1\right\rfloor}
\newcommand{\ceil}[1]{\left\lceil #1\right\rceil}

\title[A Reider-type theorem for syzygies of cyclic covers]{A Reider-type theorem for syzygies of cyclic\\ covers
of anticanonical rational surfaces}
\author{Lei Song}

\hypersetup{
  pdftitle={A Reider-type theorem for syzygies of cyclic covers of anticanonical rational surfaces},
  pdfauthor={Lei Song},
  pdfsubject={Higher syzygies of cyclic coverings},
  pdfkeywords={cyclic cover, property Np, Koszul cohomology, rational surface}
}

\begin{document}

\begin{abstract}
Let $\pi\colon X\to S$ be a smooth ramified cyclic cover of an
anticanonical rational surface, with $\kappa(X)\geq0$, and let $A$ be a divisor on $S$. We prove a Reider-type theorem for syzygies: for every integer $p\ge 0$, the adjoint divisor $K_X+\pi^*A$ satisfies Property $\mathrm{N}_p$ provided that $A\cdot C\geq p+3$ for every curve $C\subseteq S$ with $C^2\le 1$. One ingredient is a syzygy criterion for finite covers that reduces the weight-two Koszul vanishing on $X$ to weight-one vanishing for the nontrivial character modules on $S$. The most delicate case, in which $\dim|-K_S|=0$ and $K_S^2<0$, is handled by a separate geometric argument.
\end{abstract}
\maketitle

\begingroup\small\tableofcontents\endgroup

\section{Introduction}
\subsection{Main results}
Reider's theorem \cite{Reider88} gives numerical criteria for an adjoint linear
series on a smooth projective surface $X$ to be basepoint-free or very ample. Specifically, if $A$ is a nef divisor with $A^2$ sufficiently large, then the failure of $|K_X+A|$ to be free or very ample is detected by special curves $C$ on which $A$ has small degree $A\cdot C$. It is natural to ask whether the same principle governs higher syzygies.
\begin{question}\label{existence of Reider-type syzygy theorem}
Is there an analogue of Reider's criterion for syzygies of surfaces? More precisely, do the numbers $A^2$ and $A\cdot C$, as $C$ ranges over curves on $X$, control the degrees of equations defining the embedding by $|K_X+A|$ and of the relations among them?
\end{question}

Recall Property $\NP_p$ of Green and Lazarsfeld; see \cite[Section~6.1]{EL26}. Let $L$ be a very ample line bundle on a projective variety $X$, and let $X\hookrightarrow\prj{N}$ be the embedding defined by the complete linear series $|L|$. Then $L$ satisfies $\NP_0$ if $L$ is normally generated, that is, its section ring \(R(X,L)\) is generated in degree
one, and $\NP_1$ if in addition the homogeneous ideal of $X$ is generated by quadrics. More generally, for $p\ge 1$, $L$ satisfies $\NP_p$ if it satisfies $\NP_0$ and the minimal free resolution of \(R(X,L)\), viewed as a graded module over the homogeneous coordinate ring of $\prj{N}$, is linear up to the first $p$ steps. A Cartier divisor $D$ satisfies Property $\NP_p$ if the line bundle $\sshf{X}(D)$ does.

For K3 and abelian surfaces, Reider-type criteria for Property $\NP_p$ are known; see \cite{AKL19, KL19, Ito18, LSW26}. A general answer to Question \ref{existence of Reider-type syzygy theorem}, however, remains out of reach.

A smooth rational surface \(S\) is \emph{anticanonical} if \(|-K_S|\ne\varnothing\). This class includes the del Pezzo surfaces, the blow-ups of Hirzebruch surfaces at no more than eight points, and all smooth projective toric surfaces. On such a surface, basepoint-freeness and Property $\NP_p$ of an ample line bundle $L$ admit particularly simple numerical descriptions, due respectively to Harbourne \cite[Theorem III.1]{Harbourne} and to Gallego--Purnaprajna \cite[Theorem 1.3]{GP01}:
\begin{itemize}
    \item[(i)] $L$ is basepoint-free if and only if $L\cdot (-K_S)\ge 2$.
    \item[(ii)] For $p\ge 0$, $L$ satisfies Property $\NP_p$ if and only if $L\cdot (-K_S)\ge p+3$.\footnote{The necessity in (ii) requires sectional genus $g(L)>0$.}
\end{itemize}

In this paper we establish a Reider-type criterion for syzygies of cyclic covers of anticanonical rational surfaces, continuing the work of \cite{RS17, Song}.

Surfaces admitting a covering structure over a rational surface are ubiquitous in the geography and moduli of surfaces of general type. For instance, if $X$ is a Horikawa surface \cite{Horikawa76}, then the canonical system $|K_X|$ exhibits $X$, or its canonical model, as a double cover of $\prj{2}$ or of a Hirzebruch surface $\mathbb{F}_e$.

Our main result is the following.
\begin{theorem}\label{thm:main}
Let
\(X\) be a smooth projective surface with $\kappa(X)\ge 0$, and let \(\pi\colon X\to S\) be a
ramified cyclic cover of an anticanonical rational surface. Let \(A\) be a divisor on \(S\) and \(p\ge 0\) an integer, and suppose that
\begin{equation}\label{eq:curve-lower-bound}
 A\cdot C\geq p+3
 \qquad\text{for every curve }C\subseteq S \text{ with }C^2\le 1.
\end{equation}
Then \(K_X+\pi^*A\) satisfies Property \(\NP_p\).
\end{theorem}

Note that, by \cref{lem:small-self-intersection}, hypothesis \eqref{eq:curve-lower-bound} is equivalent to the a priori stronger condition $A\cdot C\ge p+3$ for \emph{every} curve $C\subseteq S$. Moreover, it implies that \(A\) is ample and basepoint-free, with \(A^2\geq(p+3)^2\); see \cref{prop:volume}.

\cref{thm:main} extends the projective normality theorem of \cite{Song}
to all \(p\geq0\). It also allows $\dim |-K_S|=0$, with the unique anticanonical divisor possibly reducible or nonreduced. Neither the covering divisor nor
\(-K_S\) need be nef.

In particular, \cref{thm:main} confirms Mukai's conjecture
\cite[Conjecture~4.2]{EL93} for smooth nonruled cyclic covers
\(\pi\colon X\to S\) of anticanonical rational surfaces,
with respect to polarizations \(\pi^*H\), where \(H\) is
an ample divisor on \(S\).
\begin{corollary}
\label{cor:Song-extension}
Let \(S\) be an anticanonical rational surface and
\(\pi\colon X\to S\) be a smooth ramified cyclic
cover, and suppose that \(X\) is nonruled. If \(H\) is an ample divisor on \(S\), then
\[
 K_X+k\pi^*H
\]
is very ample and satisfies Property \(\NP_p\) for all integers \(p\ge 0\) and
\(k\geq p+3\).\qed
\end{corollary}

By Bangere--Lacini \cite[Theorem 1.2]{BangereLacini}, if $(X, H)$ is a smooth polarized surface with $H$ basepoint-free, then \(K_X+kH\) satisfies Property $\NP_p$ for $k\ge p+3$. In \cref{cor:Song-extension}, by contrast, $H$ is not assumed to be basepoint-free, so $\pi^*H$ need not be basepoint-free either.

\subsection{Notation and outline of the proof}

We now fix the notation used throughout the paper, and then outline the proof.

Recall that a ramified cyclic cover \(\pi\colon X\to S\) of degree $r\ge 2$ is determined by a Cartier divisor \(B\) on \(S\) and a section \(\sigma\in H^0(S,\sshf{S}(rB))\); one has
\begin{equation}\label{eq:cyclic-algebra}
 \pi_*\sshf{X}\iso \bigoplus_{i=0}^{r-1}\sshf{S}(-iB).
\end{equation}

The branch divisor of $\pi$ is $\Gamma=(\sigma=0)\in|rB|$. Since $S$ is smooth, $X$ is smooth if and only if \(\Gamma\) is smooth. In that case the branch divisor need
not be connected; its irreducible components are pairwise disjoint.
The canonical bundle formula gives
\[
 K_X=\pi^*\bigl(K_S+(r-1)B\bigr).
\]
Put
\begin{equation}\label{eq:notation}
 D=K_S+A+(r-1)B,\qquad D_j=D-jB
 \quad(0\leq j<r),
\end{equation}
and write \(L=\pi^*D=K_X+\pi^*A\).

For a finite-dimensional vector space \(V\) and a graded
\(\Sym V\)-module \(M=\oplus_{i\in\mathbb{Z}} M_i\), the Koszul group \(K_{a,q}(M;V)\) is the middle cohomology
of the Koszul complex
\[
 \bigwedge^{a+1}V\otimes M_{q-1}\longrightarrow
 \bigwedge^aV\otimes M_q\longrightarrow
 \bigwedge^{a-1}V\otimes M_{q+1};
\]
see \cite[Sections~1.2--1.4]{AproduNagel}.

For a Cartier divisor \(H\) on a scheme \(Y\), the section ring is
\[
R(Y,H)=\bigoplus_{i\geq0}H^0\bigl(Y,\sshf{Y}(iH)\bigr).
\]
We abbreviate
\[
 K_{a,q}(Y,H)=
 K_{a,q}\bigl(R(Y,H);H^0(Y,\sshf{Y}(H))\bigr).
\]
In these terms, \(H\) satisfies \(\NP_p\) if it is very ample and, for every \(0\leq a\leq p,\ q\geq2\),
\[
 K_{a,q}(Y,H)=0.
\]

Throughout the paper, by the \emph{cyclic-cover setup}, we mean a smooth ramified cyclic cover $\pi: X\rightarrow S$ of an anticanonical rational surface with $\kappa(X)\ge 0$ (that is, $X$ is nonruled), together with a divisor $A$ on $S$ satisfying \eqref{eq:curve-lower-bound} for a fixed integer $p\ge 0$. The covering data $B, \Gamma$ and the divisors $D, D_j$ are as in
\eqref{eq:cyclic-algebra}, and \eqref{eq:notation}.

The proof separates into a homological reduction and a geometric argument. On the homological side, Green duality shows that the relevant Koszul groups of weight at least three vanish (\cref{prop:high-weights}), and a syzygy criterion for finite covers (\cref{thm:transfer}) reduces the remaining weight-two groups to weight-one syzygies of the nontrivial character modules on $S$
\[
 N_{D,jB}=\bigoplus_{q\geq0}
 H^0\bigl(S,\sshf{S}(D_j+qD)\bigr),\qquad 1\leq j<r.
\]
The geometric part of the proof is to verify the vanishing of $K_{a, 1}\paren{N_{D,jB}; H^0(S, \sshf{S}(D))}$ for $0\le a\le p$, which we do by applying the effective theorem of Niu--Park \cite{NiuPark} to an incomplete linear series on a smooth hyperplane section. What this requires on the surface (\cref{prop:weight-one}) is $p$-very ampleness of \(D_j\), certain cohomology vanishings for line bundles, and the specialty bound
\[
 (-K_S)\cdot D_j-h^1(S,\sshf{S}(jB))\geq p+2.
\]
The smoothness of the branch divisor determines the fractional negative part
of \(B\), and this gives the $p$-very ampleness of $D_j$ needed even when
\(B\) is not nef; see \cref{prop:rounded-adjunction}.

For regular covers, the specialty bound follows from the character decomposition and restriction to an anticanonical divisor. This is where irregularity \(q(X)\) enters. By \cref{lem:regularity} and Serre duality, \(q(X)=\sum_{j=1}^{r-1}h^1(S,\sshf{S}(K_S+jB))\), and the restriction argument needs the vanishing of these groups. Kawamata--Viehweg vanishing provides it as soon as the positive part \(P_B\) of the Zariski decomposition of \(B\) is big, that is, \(P_B^2>0\); this includes every cover of general type.

When \(P_B^2=0\), the method breaks down, and for a good reason: \(P_B\) is then nef but not big, it defines a genus-one fibration, and \(q(X)\) can be arbitrarily large (\cref{ex:large-irregularity}). One must use the fibration instead. If the anticanonical system moves, then \(S\) is a rational elliptic surface, and we argue with complete and relative kernel bundles (\cref{sec:elliptic}); if it is fixed and \(K_S^2=0\), then \(S\) is a Halphen surface of index at least two, whose unique anticanonical curve does not move (\cref{sec:halphen}). Both arguments exploit the fact that \(-K_S\) is proportional to the fiber class. Thus an anticanonical curve \(F\) is nef with \(F^2=0\), so adding copies of \(F\) to a divisor one at a time keeps it ample and does not change its degree on \(F\). At each step, the restriction sequence to \(F\) reduces the required vanishing on \(S\) to a vanishing on a genus-one curve, where it can be checked directly: on a smooth fiber in the first case, and on the fixed, possibly singular, curve \(F\) in the second.

The most delicate case is the one in which
\[
 K_S^2<0,\qquad h^0(S,\sshf{S}(-K_S))=1,\qquad P_B^2=0.
\]
This case is different in kind. Here \(S\) is obtained from a Halphen surface by blowing up \(n=-K_S^2\) points on successive anticanonical divisors (\cref{lem:fixed-geometry}). The unique anticanonical divisor \(F\) may be reducible or nonreduced, each of its components has negative square, and \(-K_S\) is no longer proportional to the fiber class. Thus neither fibration argument applies: there are no anticanonical fibers to add, and since vertical blow-ups change the anticanonical divisor by the relative canonical divisor, the statement on \(S\) does not follow formally from \cref{prop:halphen-induction} applied on the Halphen surface. Nor does the restriction argument apply, because \(P_B\) is not big.

The idea is to subtract a well-chosen curve, which supplies the missing bigness, and to add it back at the end. Let \(\beta\colon S\to Y\) be the contraction to the relatively minimal model of the genus-one fibration defined by \(P_B\), a Halphen surface, and put \(H=\beta^*(-K_Y)\), a nef integral class. We subtract from \(A\) as large a multiple of \(H\) as the bound \(p+3\) permits. The resulting
class \(J=A-\tau H\) is ample, with
\[
 J\cdot C\geq p+3\quad\text{for every }C,\qquad J^2\geq3(p+3)^2.
\]
The curve \(E\) at which equality holds is a horizontal smooth rational curve with \(E^2=-1\) or \(-2\). The Hodge index theorem controls \(E\cdot C\) for every curve \(C\), which gives the required positivity of the integral divisor \(U=A-E\).    The separation criteria of Arcara--Bertram and Langer then show that each \(D_j-E\), \(1\leq j<r\), is nef and big, and \(p\)-very ample. The curve \(E\) also supplies the missing bigness: for a suitable \(0<\lambda_j\leq1\), \(jP_B+\lambda_jE\) is big, although \(jP_B\) is not, so Kawamata--Viehweg vanishing for klt pairs applies to \(K_S+jB+E\) (\cref{lem:fixed-augmented-specialty}).

Restricting to the full anticanonical divisor \(F\), we obtain the specialty bound
\[
 F\cdot (D_j-E)-h^1(S,\sshf{S}(jB+E))\geq p+2.
\]
This yields the weight-one vanishings once \(E\) has been subtracted. To add \(E\) back, note that the complete series of \(D\) restricts surjectively to \(E\simeq\prj{1}\), where the kernel bundle splits into trivial summands and summands of degree \(-1\); the restriction sequence then recovers the required vanishings for \(D_j\).

The table below lists the results used in each geometric case.
\begin{center}
\small
\begin{tabular}{@{}p{0.43\textwidth}p{0.45\textwidth}@{}}
\toprule
Case & Result used \\
\midrule
\(\rho(S)\leq2\) & \Cref{prop:low-rank} \\
\(q(X)=0\), or \(P_B^2>0\) &
  \Cref{prop:regular-case,lem:regularity} \\
\(P_B^2=0\), \(\dim|-K_S|\geq1\) & \Cref{cor:moving-anti} \\
\(\dim|-K_S|=0\), \(K_S^2=0\) &
  \Cref{thm:nonnegative} \\
\(\dim|-K_S|=0\), \(K_S^2<0\) &
  \Cref{thm:fixed-uniform} \\
\bottomrule
\end{tabular}
\end{center}

The paper is organized as follows. \Crefrange{sec:positivity}{sec:branch} develop the numerical and homological tools. \Cref{sec:regular} treats regular covers, while \cref{sec:elliptic,sec:halphen} treat rational elliptic surfaces and Halphen surfaces of arbitrary index. The most delicate case, in which $\dim|-K_S|=0$, $K_S^2<0$ and $P_B^2=0$, is the subject of \cref{sec:fixed}. Finally, \cref{sec:completion} treats surfaces \(S\) of Picard number at most two and completes the proof of \cref{thm:main}.

\noindent\textbf{Conventions.} We work throughout over the complex numbers \(\mathbb{C}\). Varieties are integral, separated schemes of finite type over \(\mathbb{C}\). Unless otherwise stated, surfaces are smooth and projective, and covers are connected. A curve on a smooth surface is a nonzero effective divisor, possibly reducible or nonreduced; irreducible curves are taken to be reduced. A cluster on a scheme $Y$ is a zero-dimensional subscheme of $Y$.
A line bundle $L$ on $Y$ is \(k\)-very ample
if its sections separate every cluster of length
\(k+1\); that is, the natural map
\[H^0(Y, L)\rightarrow H^0(Y, L|_Z)\]
is surjective. For a globally generated line bundle $L$ on $Y$, the kernel bundle $M_L$ is the kernel of the evaluation map
\[
  0\longrightarrow M_L\longrightarrow
  H^0\bigl(Y,L\bigr)\otimes\sshf{Y}
  \longrightarrow L\longrightarrow0.
\]

For a
sheaf
\(\shf{F}\) of \(\sshf{Y}\)-modules, \(h^i(Y,\shf{F})\) denotes the
dimension of \(H^i(Y,\shf{F})\). The irregularity is $q(Y)=h^1(Y, \sshf{Y})$, and $Y$ is called regular if $q(Y)=0$. For a divisor \(D\) on \(Y\), its
associated line bundle is denoted by \(\sshf{Y}(D)\).

We write $\sim, \sim_{\mathbb{Q}}$ and $\equiv$ for linear, $\mathbb{Q}$-linear and numerical equivalence. On a rational surface, linear and numerical  equivalence coincide. We do not
distinguish between a divisor and its associated line bundle when no confusion can arise.

\vspace{5pt}

\noindent\textbf{Acknowledgements. }This project grew out of a question
on Horikawa surfaces that Purna Bangere posed to me. I would like to thank him for his encouragement and useful discussions. I also benefited greatly from many discussions with Runjian Huang. During the preparation of this article, the author was partially supported by NSFC grant No.~12471043 and by Guangdong Basic and Applied
Basic Research Foundation grant No.~2025A1515012258. 

\vspace{5pt}

\noindent\textbf{AI disclosure.} 
The author formulated the problem in its present generality and devised the overall proof strategy. A generative AI tool was used, under the author’s direction, to develop some arguments within this framework and to improve the exposition. The paper was written by the author, who takes full responsibility for the content.

\section{Positivity, adjoint separation, and anticanonical curves}\label{sec:positivity}

We begin with the numerical consequences of the hypothesis on \(A\): ampleness, a lower bound for \(A^2\), and the adjoint separation statements that will be used repeatedly below.

\begin{proposition}\label{prop:volume}
Let \(S\) be an anticanonical rational surface.
Suppose that a divisor \(A\) satisfies
\begin{equation}\label{eq:A-m}
  A\cdot C\geq m
  \qquad\text{for every irreducible curve }C\subset S,
\end{equation}
where \(m\geq2\). Then \(A\) is ample and basepoint-free,
\(H^1(S,\sshf{S}(A))=0\), and
\begin{equation}\label{eq:A-square}
  A^2\geq m^2.
\end{equation}
If \(\rho(S)=2\), the stronger inequality \(A^2\geq2m^2\) holds.
\end{proposition}

\begin{proof}
The hypothesis implies that \(A\) is nef. If \(A^2=0\), then
\(H^0(S,\sshf{S}(K_S-nA))=0\) for \(n>0\), while Riemann--Roch gives
\[
  \chi(\sshf{S}(nA))=1+\frac n2 A\cdot(-K_S)>0.
\]
Thus \(nA\) is effective, and hence \(A\cdot nA>0\) by \eqref{eq:A-m}, a contradiction.
Therefore \(A^2>0\), and by Nakai--Moishezon criterion, \(A\) is ample.

If \(\rho(S)=1\), then \(S\simeq\prj{2}\), and the assertion is immediate.
If \(\rho(S)=2\), then \(S\simeq\mathbb{F}_e\). With negative section
\(C_0\) and fiber \(f\), write \(A=aC_0+bf\), and set
\[
  c=A\cdot C_0=b-ae.
\]
Since \(a=A\cdot f\geq m\) and \(c\geq m\), one has
\[
  A^2=-ea^2+2ab=ea^2+2ac\geq(e+2)m^2\geq2m^2.
\]

Assume now that \(\rho(S)\geq3\), and put \(Q=A+mK_S\). We claim that
\(Q\) is nef. By the cone theorem,
\(\overline{\mathrm{NE}}(S)\) is generated by its \(K_S\)-nonnegative
part and its \(K_S\)-negative extremal rays. On the former,
\(Q\cdot z=A\cdot z+mK_S\cdot z\geq0\). A \(K_S\)-negative extremal
contraction has neither a point nor a curve as target: these would force \(\rho(S)=1\) or \(\rho(S)=2\), respectively. It is therefore
birational. Any irreducible contracted curve \(R\) has \(R^2<0\)
by negative definiteness of the exceptional intersection form, and
\(K_S\cdot R<0\) because the ray is \(K_S\)-negative. Adjunction gives
\[
 -2\leq 2\pa(R)-2=R^2+K_S\cdot R\leq-2.
\]
Both intersection numbers are negative integers, so
\(R^2=K_S\cdot R=-1\) and \(\pa(R)=0\). Thus \(R\) is a smooth
rational \((-1)\)-curve generating the ray. Hence
\[
  Q\cdot R=A\cdot R-m\geq0.
\]
This proves the claim. Consequently,
\[
  0\leq Q\cdot A=A^2-mA\cdot(-K_S),
\]
and \eqref{eq:A-square} follows from
\(A\cdot(-K_S)\geq m\).

Finally, we apply Harbourne's theorem for nef divisors on anticanonical rational
surfaces \cite{Harbourne}. Since
\(A\cdot(-K_S)\geq m\geq2\), it gives
\(H^1(S,\sshf{S}(A))=0\) and basepoint-freeness.
\end{proof}

\begin{lemma}
\label{lem:small-self-intersection}
Let \(S\) be an anticanonical rational surface. The additive monoid
of effective divisor classes on \(S\) is generated by the classes of
irreducible curves \(C\) with \(C^2\leq1\). Consequently, for a divisor \(A\) and a real number \(m>0\), the following conditions are equivalent:
\begin{itemize}
    \item \(A\cdot C\geq m\) for every irreducible curve \(C\) with \(C^2\leq1\).
    \item \(A\cdot C\geq m\) for every curve $C$ on \(S\).
\end{itemize}
\end{lemma}

\begin{proof}
We first prove the generation statement. If \(\rho(S)=1\), then
\(S\simeq\mathbb P^2\), and the effective monoid is generated by
the class of a line. If \(\rho(S)=2\), then
\(S\simeq\mathbb F_e\) for some $e\ge 0$, and the effective monoid is generated by
the classes of a minimal section and a fiber, whose squares are $-e$ and $0$, respectively.

Suppose that \(\rho(S)\geq3\). Fix an \(F\in|-K_S|\)
and an ample divisor \(T\). By
\cite[Lemma 4.1]{LahyaneHarbourne}, every effective divisor
class has an expression
\[
    D\sim aF+\sum_i b_iN_i,
    \qquad a,b_i\in\mathbb Z_{\geq0},
\]
where \(N_i\) are irreducible curves with \(N_i^2<0\).
If \(T\cdot D<T\cdot F\), then \(a=0\), since otherwise
intersection with \(T\) would give \(T\cdot D\geq T\cdot F\).
It therefore remains to express the class of \(F\) using
irreducible curves of square at most one.

If \(F\) is reducible or nonreduced, every irreducible component
\(C\) of \(F\) satisfies
\[
    T\cdot C<T\cdot F.
\]
Thus each component class, and hence the class of \(F\), is a
nonnegative integral combination of negative curve classes.
If \(F\) is integral with \(F^2\leq1\), its class is already
an allowed generator.

Finally, suppose that \(F\) is integral with \(d=F^2\geq2\).
Then \(F\) is nef. Since \(\rho(S)\geq3\), the rational surface
\(S\) is nonminimal and contains a smooth rational
\((-1)\)-curve \(E\). Adjunction gives \(F\cdot E=1\).
Riemann--Roch yields
\[
    \chi\bigl(\mathcal O_S(F-E)\bigr)=d-1>0.
\]
Moreover, since \(F\) is nef and $F\cdot(-2F+E)=-2d+1<0$, by Serre duality,
\[
    h^2\bigl(S,\mathcal O_S(F-E)\bigr)
    =h^0\bigl(S,\mathcal O_S(-2F+E)\bigr)=0.
\]
Consequently, \(F-E\) is linearly equivalent to an effective
divisor. Since \(T\cdot(F-E)<T\cdot F\), its class is a
nonnegative integral combination of negative curve classes. So the class $F=(F-E)+E$ has the required expression. This proves the generation statement.

The second condition clearly implies the first. Conversely,
assume the first condition, and let \(C\) be any curve on \(S\).
By the generation statement,
\[
    C\sim\sum_{i=1}^{k}n_iC_i,
    \qquad n_i\in\mathbb Z_{>0},
\]
where the \(C_i\) are irreducible curves with $C_i^2\leq1$. Therefore
\[
    A\cdot C
    =\sum_{i=1}^{k}n_i A\cdot C_i
    \geq m\sum_{i=1}^{k}n_i
    \geq m.
\]
Hence the second condition holds.
\end{proof}
By \cref{lem:small-self-intersection}, the hypothesis \eqref{eq:curve-lower-bound} amounts to the bound \(A\cdot C\geq p+3\) for every irreducible curve \(C\subset S\), and hence for every nonzero curve \(C\). This is the form in which we shall use it.

\begin{lemma}\label{lem:adjoint-separation}
Let \(p\geq0\), and let \(S\) and \(A\) be as in
\cref{prop:volume}, with \(A\cdot C\geq p+3\) for every irreducible
curve \(C\subset S\).
\begin{enumerate}
\item \(K_S+A\) is \(p\)-very ample.
\item If \(f\colon Y\to S\) is a finite map of degree at least two and \(Y\) is
smooth, then \(K_Y+f^*A\) is \((p+1)\)-very ample.
\end{enumerate}
\end{lemma}

\begin{proof}
We use the following form of the Reider criterion of Arcara--Bertram
\cite[Application 2.2 and the ensuing Reider theorem]{ArcaraBertram}.
If \(H\) is ample on a smooth surface, then \(K+H\) separates every
cluster of length at most \(d\), provided
\begin{equation}\label{eq:AB-Reider}
 H^2>(d+1)^2,\qquad
 H\cdot C>C^2+d
\end{equation}
for every effective divisor \(C\) with \(C^2\leq0\).

For (i), take \(H=A\) and \(d=p+1\). By
\cref{prop:volume}, \(A^2\geq (p+3)^2>(p+2)^2=(d+1)^2\).
For an effective \(C=\sum n_iC_i\) with \(C^2\leq0\),
\[
 A\cdot C\geq (p+3)\sum n_i\geq p+3>d+C^2,
\]
so \eqref{eq:AB-Reider} applies.

For (ii), take \(H=f^*A\) and \(d=p+2\). Then
\[
  (f^*A)^2=(\deg f)A^2\geq2(p+3)^2>(p+3)^2=(d+1)^2.
\]
If \(T=\sum n_iT_i\) is an effective divisor on \(Y\) with \(T^2\leq0\),
then
\[
 f^*A\cdot T=A\cdot f_*T\geq (p+3)\sum n_i\deg(T_i/f(T_i))
 \geq p+3>d+T^2.
\]
The same criterion proves (ii).
\end{proof}

The Arcara--Bertram criterion covers the integral adjoint divisors of \cref{lem:adjoint-separation}. The rounded adjoints of
\cref{prop:rounded-adjunction} require a logarithmic
version of Reider's criterion for nef \(\mathbb{Q}\)-divisors, due to Langer \cite{Langer}.

We will use the following bound for
Langer's cluster
invariant on a smooth surface. For a cluster \(Z\subset S\), set
\begin{equation}\label{eq:Langer-cluster-bound}
 \delta_Z:=4\operatorname{length}(Z).
\end{equation}
If separation fails, a minimal
subcluster \(Z'\subseteq Z\) on which it fails is Gorenstein.
Langer's invariant for \(Z'\) is defined locally: at each \(x\in\operatorname{Supp}(Z')\), the rank-two reflexive module
\(\mathcal E_{Z',x}\) arises from a class in
\(
\operatorname{Ext}^1_{\sshf{S,x}}
       (\is{Z',x},\omega_{S,x}).
\)
The associated local invariants are then summed over \(\operatorname{Supp}(Z')\).
On a smooth surface the resulting invariant equals
\(4\operatorname{length}(Z')\leq\delta_Z\); see
\cite[(0.3.2) and Definitions 1.1.1--1.1.2]{Langer}

\begin{lemma}[Langer]
\label{lem:Langer-rounded}
Let \(S\) be a smooth projective surface, let \(M\) be a nef
\(\mathbb{Q}\)-divisor whose fractional support has simple normal crossings,
and let \(Z\subset S\) be a cluster. Suppose that positive rational
numbers \(\alpha,\beta\) satisfy
\begin{equation}\label{eq:Langer-parameters}
 \alpha\geq\delta_Z,
 \qquad 4\beta\left(1-\frac{\beta}{\alpha}\right)\geq\delta_Z,
\end{equation}
and that
\begin{equation}\label{eq:Langer-numerics}
 M^2>\alpha,
 \qquad M\cdot C\geq\beta
 \quad\text{for every irreducible curve }C\subset S.
\end{equation}
Then the restriction map
\[
 H^0\bigl(S,K_S+\ceil M\bigr)\longrightarrow
 H^0\bigl(Z,(K_S+\ceil M)|_Z\bigr)
\]
is surjective.
\end{lemma}

\begin{proof}
If separation fails, choose a minimal failing subcluster
\(Z'\subseteq Z\). Its Langer invariant is at most \(\delta_Z\),
so \eqref{eq:Langer-parameters} also holds for that invariant.
The smooth log-terminal case of \cite[Theorem 0.2]{Langer}
contradicts failure of separation on \(Z'\). In that theorem only
the square inequality in \eqref{eq:Langer-numerics} is strict;
the curve inequality and both parameter inequalities are weak.
The simple-normal-crossing hypothesis ensures that the boundary
\(\ceil M-M\) has the required singularities.
\end{proof}

\begin{lemma}[Anticanonical subcurves]\label{lem:anti-subcurves}
Let \(S\) be a smooth rational surface and let \(F\in|-K_S|\).
Then
\begin{itemize}
    \item[(i)] \(H^0(F,\sshf{F})=\mathbb{C}\), \(\pa(F)=1\), and
\(\omega_F\simeq\sshf{F}\).
    \item[(ii)] Every proper positive subcycle \(Y<F\)
satisfies \(H^1(Y,\sshf{Y})=0\), and hence \(\pa(Y)\leq0\).
    \item[(iii)] If \(h^0(S,\sshf{S}(-K_S))=1\), then \(K_S^2\leq0\).
\end{itemize}
\end{lemma}
\begin{proof}
The anticanonical sequence
\[
 0\longrightarrow\sshf{S}(K_S)\longrightarrow\sshf{S}
 \longrightarrow\sshf{F}\longrightarrow0
\]
and rationality give the first assertions, together with adjunction.
For \(0<Y<F\), the divisor \(K_S+Y\sim-(F-Y)\) is not effective.
The sequence
\[
 0\longrightarrow\sshf{S}(K_S)\longrightarrow\sshf{S}(K_S+Y)
 \longrightarrow\omega_Y\longrightarrow0
\]
therefore gives \(H^0(Y,\omega_Y)=0\). Duality on the Cartier curve
\(Y\) yields the assertion about \(H^1\). Finally, Riemann--Roch and
\(H^0(S,\sshf{S}(2K_S))=0\) give
\(h^0(S,\sshf{S}(-K_S))=1+K_S^2
+h^1(S,\sshf{S}(-K_S))\), proving the last assertion.
\end{proof}

\section{Homological reductions}\label{sec:homological}

We first dispose of the Koszul groups of weight at
least three, and then prove a result that reduces the weight-two Koszul vanishing to weight-one vanishing of the nontrivial character summands. These reductions will allow us to treat the geometry of the covering data separately.

\begin{proposition}[Vanishing in higher weights]\label{prop:high-weights}
In the cyclic-cover setup, we have
\[
  K_{a,q}(X,L)=0\qquad \text{ for all } 0\leq a\leq p,\ q\geq3.
\]
\end{proposition}

\begin{proof}
Put \(P=\pi^*A\), so \(L=K_X+P\). By
\cref{lem:adjoint-separation}\textup{(ii)}, \(L\) is very ample.
By Kodaira vanishing, \(H^1(X,\sshf{X}(jL))=0\) for \(j\geq1\),
since
\[
 jL=K_X+(j-1)L+P.
\]
Set \(R=h^0(X,\sshf{X}(L))-1\). Green duality \cite{Green} gives, for $a\ge 0$ and $q\ge 3$,
\begin{equation*}
 K_{a,q}(X,L)\simeq
 K_{R-a-2,\,3-q}(X,K_X;L)^\vee.
\end{equation*}
For \(q\geq4\), the group on the right has negative weight, and vanishes because
\(K_X+(3-q)L=-P-(q-4)L\) has no nonzero sections.
For \(q=3\), Green's vanishing theorem
\cite[Theorem~(3.a.1)]{Green} gives
\(K_{s,0}(X,K_X;L)=0\) for \(s\geq p_g(X)\). For \(0\leq a\leq p\), it therefore suffices to prove
\begin{equation}\label{eq:h0-pg}
 h^0(X,\sshf{X}(L))-p_g(X)\geq p+3.
\end{equation}

Let \(g(P)=1+(K_X+P)\cdot P/2\) be the sectional genus of the
ample, globally generated divisor \(P\). The Riemann--Roch theorem and Kodaira
vanishing give
\[
 h^0(X,\sshf{X}(L))-p_g(X)=g(P)-q(X).
\]
Because \(X\) is nonruled, it is not a scroll over a curve of genus at
least two. The Fukuma--Ishihara inequality
\cite[Proposition 4.1]{FukumaIshihara} therefore gives
\(g(P)\geq2q(X)-1\). Moreover, \(K_X\) is pseudoeffective, so
\(K_X\cdot P\geq0\). Hence
\begin{align*}
 h^0(X,\sshf{X}(L))-p_g(X)
 &=g(P)-q(X)\\
 &\geq\frac{g(P)-1}{2}
 =\frac{(K_X+P)\cdot P}{4}
 \geq\frac{P^2}{4}
 =\frac{rA^2}{4}.
\end{align*}
By \cref{prop:volume}, the last quantity is at least
\((p+3)^2/2\geq p+3\), proving \eqref{eq:h0-pg}. This completes the proof.
\end{proof}

It remains to treat the weight-two Koszul groups. The following
theorem reduces their vanishing on a finite cover to syzygy
conditions on the base.

\begin{theorem}[Syzygy criterion for finite covers]\label{thm:transfer}
Let \(f\colon Y\to Z\) be finite and surjective between
projective varieties, with
\[
 f_*\sshf{Y}\iso \sshf{Z}\oplus\shf{E}
\]
as \(\sshf{Z}\)-modules. Let \(D\) be ample, put \(L=f^*D\), and set
\[
 V=H^0\bigl(Z,\sshf{Z}(D)\bigr),\qquad
 N_D(\shf{E})=
 \bigoplus_{q\geq0}H^0\!\left(Z,\shf{E}\otimes\sshf{Z}((q+1)D)\right),
\]
where the displayed summand has degree \(q\). Fix \(p\geq0\), and
assume that:
\begin{enumerate}
\item \(D\) satisfies \(\NP_p\);
\item \(K_{a,1}(N_D(\shf{E});V)=0\) for \(0\leq a\leq p\);
\item \(K_{a,q}(Y,L)=0\) for \(0\leq a\leq p\) and \(q\geq3\).
\end{enumerate}
Then \(L\) is very ample and satisfies \(\NP_p\). If \(\shf{E}\) is a
finite direct sum, condition \textup{(ii)} can be checked summand by
summand.
\end{theorem}

\begin{proof}
Write \(U=H^0(Y,\sshf{Y}(L))\), and choose a splitting
\(U=V\oplus W\). In internal degree \(d\), the Koszul complex has terms
\[
 C_n(d)=\bigwedge^nU\otimes R(Y,L)_{d-n}
 \iso \bigoplus_{a+b=n}\bigwedge^aV\otimes\bigwedge^bW
       \otimes R(Y,L)_{d-n}.
\]
The differential is the sum of multiplication by the \(V\)-sections
and multiplication by the \(W\)-sections. The former preserves the
number of \(W\)-factors, while the latter decreases it by one. Thus
the increasing filtration
\[
 \mathcal F_bC_n(d)=
 \bigoplus_{0\leq j\leq\min\{b,n\}}\bigwedge^{n-j}V\otimes\bigwedge^jW
       \otimes R(Y,L)_{d-n}
\]
is preserved by the differential. On the associated graded complex,
only the \(V\)-differential remains. For fixed \(b\), this is the
Koszul complex on \(V\) in internal degree \(d-b\), tensored with
\(\bigwedge^bW\) and shifted in homological degree by \(b\).
Taking its homology gives the first page
\begin{equation}\label{eq:change-rings}
 E^1_{a,b}(d)=
 K_{a,d-a-b}(R(Y,L);V)\otimes\bigwedge^bW
 \ \Longrightarrow\
 K_{a+b,d-a-b}(Y,L).
\end{equation}
The filtration is finite, so the spectral sequence converges to the
homology of \(C_\bullet(d)\); its limiting terms are the associated
graded pieces of that homology. The projection formula gives, as graded
\(\Sym V\)-modules,
\begin{equation}\label{eq:ring-decomp}
 R(Y,L)\iso R(Z,D)\oplus M_D(\shf{E}),
 \qquad M_D(\shf{E})=N_D(\shf{E})(-1).
\end{equation}
Here \(H^0(Z,\shf{E})=0\): taking global sections of the split unit
summand gives
\[
 H^0(Y,\sshf{Y})\iso H^0(Z,\sshf{Z})\oplus H^0(Z,\shf{E})\iso \mathbb{C}.
\]

For total homological degree \(n=a+b\leq p\) and internal degree \(n+2\),
\eqref{eq:change-rings} reads
\[
 E^1_{a,b}(n+2)=
 K_{a,2}(R(Y,L);V)\otimes\bigwedge^bW.
\]
The summand \(K_{a,2}(R(Z,D);V)\) vanishes by \textup{(i)}, while
\eqref{eq:ring-decomp} identifies the other summand with
\[
 K_{a,2}(M_D(\shf{E});V)
 =K_{a,1}(N_D(\shf{E});V),
\]
which vanishes by (ii). Hence \(E^1_{n-b,b}(n+2)=0\) for every
\(0\leq b\leq n\). Each associated graded piece of
\(K_{n,2}(Y,L)\) is a subquotient of one of these terms and is
therefore zero. Since the filtration is finite, \(K_{n,2}(Y,L)=0\)
for \(0\leq n\leq p\). Hypothesis (iii) removes all higher weights.

The case \(n=0\), together with (iii), shows that the section ring $R(Y, L)$ is
generated in degree one. Moreover, \(L=f^*D\) is basepoint-free because
\(D\) is very ample. It follows from
\cite[pp.~38--39]{Mumford} that $L$ is very ample, so the asserted Koszul-cohomology vanishing is Property
\(\NP_p\). The last assertion follows because Koszul cohomology commutes with finite direct sums.
\end{proof}

The next lemma gives a sufficient condition for the vanishing
of the weight-one Koszul groups of an incomplete twisted section
module on a smooth projective curve.

\begin{lemma}\label{lem:incomplete-weight-one}
Let \(C\) be a smooth projective curve, let \(\shf{B},\shf{L}\)
be line bundles, and let \(W\subseteq H^0(C,\shf{B})\).
Suppose that \(W\) separates every cluster of length \(p+1\),
where \(p\geq0\), and that
\begin{equation}\label{eq:incomplete-specialty}
 h^1(C,\shf{L}\otimes\shf{B}^{-1})
 \leq\dim W-p-2.
\end{equation}
Form the graded section module with degree-zero part \(W\) and
degree-\(q\) part \(H^0(C,\shf{B}\otimes\shf{L}^q)\) for
\(q\geq1\), over \(\Sym H^0(C,\shf{L})\).
Its groups \(K_{a,1}\) vanish for \(0\leq a\leq p\).
\end{lemma}

\begin{proof}
It suffices to prove the assertion for \(a=p\), since the
hypotheses remain valid for every smaller index. On the symmetric product
\(C_{p+1}\), write
\(\shf{E}_{p+1,\shf{B}}\) for the tautological bundle and
\(\shf{N}_{p+1,\shf{B}}=\det\shf{E}_{p+1,\shf{B}}\).
The separation hypothesis gives the kernel bundles
\[
 \shf{M}_W=\ker(W\otimes\sshf{C_{p+1}}\longrightarrow
                         \shf{E}_{p+1,\shf{B}}),\qquad
 \shf{M}=\ker(H^0(C,\shf{B})\otimes\sshf{C_{p+1}}\longrightarrow
                         \shf{E}_{p+1,\shf{B}}).
\]
Put \(U=H^0(C,\shf{B})/W\), \(n=\rk{\shf{M}_W}\), and
\(k=n-1=\dim W-p-2\). By the snake lemma, there is an exact sequence
\[
 0\longrightarrow\shf{M}_W\longrightarrow\shf{M}
 \longrightarrow U\otimes\sshf{C_{p+1}}\longrightarrow0.
\]
For \(\shf{T}=\omega_C\otimes\shf{B}\otimes\shf{L}^{-1}\),
condition \eqref{eq:incomplete-specialty} amounts to \(h^0(C,\shf{T})\leq k\).
Niu--Park's symmetric-product vanishing
\cite[Theorem 3.5]{NiuPark} therefore gives
\begin{equation}\label{eq:incomplete-resolution-vanishing}
 H^{p-i}\!\left(C_{p+1},
 \bigwedge^{k-i}\shf{M}\otimes\shf{N}_{p+1,\shf{T}}\right)=0
 \quad\quad \text{for all } 0\leq i\leq\min\{p,k\}.
\end{equation}
Indeed, the bound in that theorem is
\(h^0(\shf{T})\leq p+(k-i)-(p-i)=k\).

The exterior-power resolution
\[
 0\to\bigwedge^k\shf{M}_W\to\bigwedge^k\shf{M}
 \to\bigwedge^{k-1}\shf{M}\otimes U
 \to\cdots\to\Sym^kU\otimes\sshf{C_{p+1}}\to0
\]
and \eqref{eq:incomplete-resolution-vanishing} imply
\[
 H^p\!\left(C_{p+1},\bigwedge^k\shf{M}_W
                  \otimes\shf{N}_{p+1,\shf{T}}\right)=0.
\]
Since \(\det\shf{M}_W=\shf{N}_{p+1,\shf{B}}^{-1}\)
up to a constant one-dimensional factor, Serre duality identifies
this vanishing with
\[
 H^1(C_{p+1},\shf{M}_W\otimes\shf{N}_{p+1,\shf{L}})=0.
\]
Here we used
\(\omega_{C_{p+1}}\otimes\shf{N}_{p+1,\shf{B}}
\otimes\shf{N}_{p+1,\shf{L}}^{-1}
=\shf{N}_{p+1,\shf{T}}\).
Consequently
\[
 W\otimes H^0(C_{p+1},\shf{N}_{p+1,\shf{L}})
 \longrightarrow H^0(C_{p+1},
 \shf{E}_{p+1,\shf{B}}\otimes\shf{N}_{p+1,\shf{L}})
\]
is surjective. The symmetric-product description of Koszul cycles
used in \cite[Section 1]{NiuPark} identifies its cokernel with
\(K_{p,1}\) for the stated module. Only the source of this map involves degree zero.
Replacing \(H^0(C,\shf{B})\) by \(W\) therefore gives
the Koszul complex of the incomplete section module.
\end{proof}

Using the preceding lemma, we obtain a sufficient condition for
the vanishing of the weight-one Koszul groups of a twisted
section module on an anticanonical rational surface.

\begin{proposition}
\label{prop:weight-one}
Let \(S\) be an anticanonical rational surface, let \(p\geq0\), let \(D\)
be very ample, and let \(T\) be Cartier. Put \(G=D-T\), assume that
\(G\) is nef and big, and suppose
\begin{align}
 &G\ \text{is \(p\)-very ample},\notag\\
 &D\cdot T>0,\qquad
 H^0(S,\sshf{S}(-T))=H^0(S,\sshf{S}(K_S-T))=0,
 \label{eq:weight-hyp2}\\
 &H^1(S,\sshf{S}(qD-T))=0\qquad(q\geq1),
 \label{eq:weight-hyp3}\\
 &(-K_S)\cdot G-h^1(S,\sshf{S}(T))\geq p+2.\notag
\end{align}
For
\[
 N_{D,T}=\bigoplus_{q\geq0}
 H^0\bigl(S,\sshf{S}((q+1)D-T)\bigr),
\]
one then has
\[
 K_{a,1}\bigl(N_{D,T};H^0(S,\sshf{S}(D))\bigr)=0
 \qquad(0\leq a\leq p).
\]
\end{proposition}

\begin{proof}
Choose a smooth \(C\in|D|\), let \(s\in H^0(S, \sshf{S}(D))\) be a defining section, and put
\[
 \shf{L}=\sshf{C}(D),\qquad \shf{B}=\sshf{C}(G),\qquad
 W=\operatorname{im}\bigl(H^0(S,\sshf{S}(G))
 \longrightarrow H^0(C,\shf{B})\bigr).
\]
Since \(H^0(S,\sshf{S}(-T))=0\), one has
\(W\simeq H^0(S,\sshf{S}(G))\).
Restriction and \eqref{eq:weight-hyp3} give
\begin{equation}\label{eq:hyperplane-module}
 \bigl(N_{D,T}/sN_{D,T}\bigr)_q=
 \begin{cases}
 W,&q=0,\\
 H^0(C,\shf{B}\otimes\shf{L}^{q}),&q\geq1.
 \end{cases}
\end{equation}
The corresponding complete spaces of sections vanish in negative
degrees, since their degrees on \(C\) are at most \(-D\cdot T<0\).
The element \(s\) is a nonzerodivisor on \(N_{D,T}\); reducing a
minimal resolution modulo \(s\) therefore preserves all Betti
numbers. Rationality gives
\(H^0(S,\sshf{S}(D))/\mathbb{C} s=H^0(C,\shf{L})\).

The subspace \(W\) separates every length-\((p+1)\) cluster of \(C\)
because \(G\) is \(p\)-very ample on \(S\). Consider the short exact sequence
\[\ses{\sshf{S}}{\sshf{S}(D)}{\sshf{C}(D)}.\]
Tensoring it with $\sshf{S}(-G)$ and taking cohomology, we get the exact sequence
\[
\begin{gathered}
\cdots\rightarrow H^1(S, \sshf{S}(-G))\rightarrow H^1(S, \sshf{S}(T))
\rightarrow H^1(C, \shf{L}\otimes\shf{B}^{-1})\\
\rightarrow H^2(S, \sshf{S}(-G))\rightarrow H^2(S, \sshf{S}(T))\rightarrow \cdots
\end{gathered}
\]
Since \(G\) is nef and big, Kawamata--Viehweg vanishing and Serre duality give $H^1(S, \sshf{S}(-G))=0$. And Serre duality and \eqref{eq:weight-hyp2} give $H^2(S, \sshf{S}(T))=0$. It follows that
\[
 h^1(C,\shf{L}\otimes\shf{B}^{-1})
 =h^0(S,\sshf{S}(K_S+G))+h^1(S,\sshf{S}(T)).
\]
By the last hypothesis, \((-K_S)\cdot G\geq p+2\geq2\), so Harbourne's theorem gives \(H^1(S,\sshf{S}(G))=0\). Moreover, \(H^2(S,\sshf{S}(G))\simeq H^0(S,\sshf{S}(K_S-G))^\vee=0\); and \(H^i(S,\sshf{S}(K_S+G))=0\) for \(i>0\) by Kawamata--Viehweg vanishing. Hence \(h^0=\chi\) for both \(\sshf{S}(G)\) and \(\sshf{S}(K_S+G)\), and we obtain that
\begin{eqnarray*}
 &&\dim W-h^1(C,\shf{L}\otimes\shf{B}^{-1})\\
 &=&h^0(S,\sshf{S}(G))-h^0(S,\sshf{S}(K_S+G))
   -h^1(S,\sshf{S}(T))\\
 &=& \chi(\sshf{S}(G))-\chi(\sshf{S}(K_S+G))-h^1(S,\sshf{S}(T))\\
 &=&(-K_S)\cdot G-h^1(S,\sshf{S}(T))\geq p+2,
\end{eqnarray*}
where the last inequality is by the assumption.

We finish the proof by applying \cref{lem:incomplete-weight-one} to \eqref{eq:hyperplane-module}.
\end{proof}

\medskip\noindent\textit{Remark.}
Under the hypotheses of \cref{prop:weight-one}, 
\[
h^1(S,\sshf{S}(T))
=h^0(S,\sshf{S}(T))-\chi(S,\sshf{S}(T)).
\]
For this reason, we call $h^1(S,\sshf{S}(T))$ the \textit{specialty} of $\sshf{S}(T)$, and 
refer to the last inequality in the proposition as the
\emph{specialty bound}.

\begin{corollary}\label{cor:cyclic-reduction}
In the cyclic-cover setup,
suppose that \(D\) satisfies \(\NP_p\) and that
\cref{prop:weight-one} applies to \((D,T=jB)\) for every
\(1\leq j<r\). Then \(L=\pi^*D\) satisfies \(\NP_p\).
\end{corollary}

\begin{proof}
The trace-zero part is
\(\bigoplus_{j=1}^{r-1}\sshf{S}(-jB)\), so
\cref{prop:weight-one} gives condition (ii) of \cref{thm:transfer}
character by character. \Cref{prop:high-weights} gives condition (iii),
and the result follows.
\end{proof}

We close this section with two consequences of the kernel-bundle description of Koszul cohomology, which will be used in \cref{sec:elliptic,sec:halphen}. Let \(D\) be a globally generated divisor on a projective variety \(Y\), put \(V=H^0(Y,\sshf{Y}(D))\), and let \(M_D\) be the kernel bundle of \(D\). Taking exterior powers of the evaluation sequence gives exact sequences
\begin{equation}\label{eq:koszul-wedge}
 0\longrightarrow\bigwedge^{w}M_D\longrightarrow\bigwedge^{w}V\otimes\sshf{Y}
 \longrightarrow\bigwedge^{w-1}M_D\otimes\sshf{Y}(D)\longrightarrow0
 \qquad(w\geq1).
\end{equation}
For a Cartier divisor \(G\) on \(Y\), let \(N=\bigoplus_{q\geq0}H^0(Y,\sshf{Y}(G+qD))\), graded so that \(H^0(Y,\sshf{Y}(G+qD))\) has degree \(q\). Twisting \eqref{eq:koszul-wedge} by \(G+qD\) and taking cohomology identifies 
\[K_{w-1,q+1}(N;V)\subseteq H^1\bigl(Y,\bigwedge^wM_D\otimes\sshf{Y}(G+qD)\bigr)\]
as a subspace, with equality if \(H^1(Y,\sshf{Y}(G+qD))=0\); see \cite[Section~2.1]{AproduNagel}.

\begin{lemma}\label{lem:kernel-criterion}
In this situation, let \(p\geq0\), and suppose that
\begin{equation*}
H^1\!\left(Y,\bigwedge^{i+1}M_D\otimes\sshf{Y}(G+qD)\right)=0
 \qquad(0\leq i\leq p,\ q\geq0).
\end{equation*}
Then \(K_{i,j}(N;V)=0\) for \(0\leq i\leq p\) and \(j\geq1\); that is, the minimal free resolution of \(N\) is linear in the first \(p\) steps. In particular, if \(Y=S\) and \(G=D_j\) in the cyclic-cover setup, then \(N=N_{D,jB}\), and condition \textup{(ii)} of \cref{thm:transfer} holds for the summand \(\sshf{S}(-jB)\).\qed
\end{lemma}

\begin{lemma}\label{lem:weight-one-regularity}
Let \(S\) be an anticanonical rational surface, let \(m\geq0\), and suppose that \((D,T)\) satisfies the hypotheses of \cref{prop:weight-one} with \(m\) in place of \(p\). Put \(G=D-T\), and assume in addition that \(H^0(S,\sshf{S}(K_S+T))=0\). Then
\begin{equation}\label{eq:weight-one-regularity}
 H^1\!\left(S,\bigwedge^wM_D\otimes\sshf{S}(G+qD)\right)=0
 \qquad(0\leq w\leq m+1,\ q\geq0).
\end{equation}
\end{lemma}

\begin{proof}
Put \(V=H^0(S,\sshf{S}(D))\) and \(N=\bigoplus_{q\geq0}H^0(S,\sshf{S}(G+qD))\); this is the module \(N_{D,T}\) of \cref{prop:weight-one}. Since \(G+qD=(q+1)D-T\), hypothesis \eqref{eq:weight-hyp3} says that \(H^1(S,\sshf{S}(G+qD))=0\) for every \(q\geq0\). This is \eqref{eq:weight-one-regularity} for \(w=0\), and for \(w\geq1\) it identifies the group in \eqref{eq:weight-one-regularity} with \(K_{w-1,q+1}(N;V)\). When \(q=0\), these groups vanish for \(w\leq m+1\) by \cref{prop:weight-one}.

When \(q\geq1\), we use Castelnuovo--Mumford regularity with respect to the very ample divisor \(D\). The line bundle \(\sshf{S}(G)\) is \(1\)-regular, since
\[
 H^1(S,\sshf{S}(G))=0,\qquad
 H^2(S,\sshf{S}(G-D))=H^2(S,\sshf{S}(-T))\simeq H^0(S,\sshf{S}(K_S+T))^\vee=0.
\]
Moreover, \(H^0(S,\sshf{S}(G-kD))=0\) for \(k\geq1\); otherwise \(-T=(G-kD)+(k-1)D\) would be effective, contrary to \eqref{eq:weight-hyp2}. Thus \(N\) is the full module of global sections of the \(1\)-regular sheaf \(\sshf{S}(G)\), for the embedding of \(S\) defined by \(|D|\). Hence \(K_{a,j}(N;V)=0\) for all \(a\) and all \(j\geq2\), which proves \eqref{eq:weight-one-regularity} for \(q\geq1\).
\end{proof}

\section{The covering data and the character adjoints}\label{sec:branch}

We now use the smoothness of the branch divisor to control
the covering divisor. The key point is that every curve on which \(B\) is negative is a branch component, disjoint from the others and occurring with a prescribed coefficient. This determines the fractional negative part of \(B\), and it gives positivity of the character adjoints \(D_j\) even when
\(B\) is not nef.

\begin{proposition}\label{prop:fractional-negative}
Suppose that \(\Gamma\in|rB|\) is smooth. Let
\(E_1,\dots,E_s\) be all irreducible curves with
\(B\cdot E_\alpha<0\), and write \(B\cdot E_\alpha=-n_\alpha\).
Then the \(E_\alpha\) are pairwise disjoint smooth components of
\(\Gamma\), and
\[
  E_\alpha^2=-rn_\alpha.
\]
Moreover,
\begin{equation}\label{eq:zariski-B}
 B=P_B+E,\qquad
 E=\frac1r\sum_\alpha E_\alpha,\qquad
 P_B=B-E
\end{equation}
is the Zariski decomposition of \(B\). In particular, \(P_B\) is nef,
\(P_B\cdot E_\alpha=0\), and
\begin{equation}\label{eq:round-jB}
  jB=\ceil{jP_B}\qquad(1\leq j<r).
\end{equation}
In the cyclic-cover setup, \(P_B\not\equiv0\).
If \(\dim|-K_S|\geq1\), all \(E_\alpha\) are rational.
\end{proposition}

\begin{proof}
If \(B\cdot C<0\), then \(\Gamma\cdot C<0\), and hence \(C\) is a
component of \(\Gamma\). Smoothness implies that \(\Gamma\) is reduced
and that its distinct components are disjoint. Writing
\(\Gamma=C+\Gamma'\), one obtains
\[
 C^2=\Gamma\cdot C=rB\cdot C.
\]
Adjunction gives
\[
 (-K_S)\cdot E_\alpha=2-2g(E_\alpha)-rn_\alpha.
\]
If \(g(E_\alpha)>0\), this number is negative, so every
anticanonical divisor contains \(E_\alpha\).
By \cref{lem:anti-subcurves}, it cannot be a proper subcycle of an
anticanonical divisor. Thus \(F=E_\alpha\) is the unique
anticanonical divisor, and adjunction gives \(g(E_\alpha)=1\).
This proves the assertions about the genus.

Set \(P_B=B-\frac1r\sum E_\alpha\). It is orthogonal to every
\(E_\alpha\). If \(C\neq E_\alpha\) for all \(\alpha\), then
\[
 rP_B\cdot C=(\Gamma-\textstyle\sum E_\alpha)\cdot C\geq0.
\]
Thus \(P_B\) is nef. The intersection matrix of the \(E_\alpha\) is
diagonal and negative definite, proving \eqref{eq:zariski-B}.
Formula \eqref{eq:round-jB} follows because \(0<j/r<1\).

Finally, in the cyclic-cover setup,
\cref{lem:covering-Iitaka} gives \(\kappa(B)\geq1\). If
\(P_B\equiv0\), the effective divisor
\(\Gamma-\sum_\alpha E_\alpha\sim rP_B\) has zero intersection
with an ample divisor, so it must be zero. Consequently
\(P_B\sim_{\mathbb Q}0\). Since Zariski decomposition preserves
Iitaka dimension, this would give
\(\kappa(B)=\kappa(P_B)=0\), a contradiction.
\end{proof}

The following example shows that in \cref{prop:fractional-negative} the hypothesis that
\(\dim|-K_S|\geq1\) is necessary to ensure all $E_{\alpha}$ are rational. Without this assumption, a nonrational \(E_\alpha\) can occur; it is then the
unique anticanonical divisor, and a smooth elliptic curve.

\begin{example}
\label[example]{rem:elliptic-negative-example}
Let \(S\) be the blow-up of \(\prj{2}\) at the fifteen transverse
intersection points of a smooth cubic and a smooth quintic. Write
\(H\) for the pullback of a line and \(E_i\) for the exceptional
curves. The strict transforms
\[
 F=3H-\sum_{i=1}^{15}E_i=-K_S,\qquad
 G=5H-\sum_{i=1}^{15}E_i
\]
are smooth and disjoint. B\'ezout's theorem shows that \(F\) is the
unique anticanonical member, and \(F^2=-6\).
Set \(B=4H-\sum^{15}_{i=1} E_i\). Then \(2B\sim F+G\), so the smooth divisor
\(F+G\) defines a smooth double cover, and
\[
 B\cdot F=-3,\qquad P_B=\tfrac12G,\qquad
 P_B^2=\tfrac52,\qquad K_S+B=H.
\]
The cover is of general type, because \(K_X=\pi^*H\) is big.
\end{example}

We record a restriction on the Iitaka dimension of the divisor $B$.
\begin{lemma}
\label{lem:covering-Iitaka}
In the cyclic-cover setup, \(\kappa(B)\geq1\).
\end{lemma}

\begin{proof}
Since \(\Gamma\in|rB|\) is nonzero and effective,
\(\kappa(B)\geq0\). Suppose to the contrary that \(\kappa(B)=0\). Then
\(\dim |n\Gamma|=0\) for every integer \(n\geq1\).

Since \(\kappa(X)\geq0\), we can find \(k\geq1\)
such that \(H^0(X,\sshf{X}(kK_X))\neq0\). The canonical bundle
formula and the character decomposition give
\[
 H^0(X,\sshf{X}(kK_X))=
 \bigoplus_{i=0}^{r-1}H^0\!\left(S,
 \sshf{S}\bigl(kK_S+(k(r-1)-i)B\bigr)\right).
\]
A nonzero summand therefore supplies an effective divisor
\[
 Q\sim kK_S+dB,\qquad d=k(r-1)-i,
 \qquad 0\leq i<r.
\]
Fix \(F\in|-K_S|\) and an ample divisor \(H\) on
\(S\). Since \(F\cdot H>0\), \(B\cdot H>0\), and
\(Q\cdot H=-kF\cdot H+dB\cdot H\geq0\), we have
\(0<d\leq k(r-1)<rk\). Moreover,
\[
 rkF+rQ\sim d\Gamma.
\]
The assumed uniqueness of \(|d\Gamma|\) makes this an equality
of effective divisors. Let \(C\) be an irreducible component of
\(F\). It must occur in \(\Gamma\), since it occurs on the left.
Its coefficient on the left is at least \(rk\), whereas its
coefficient on the right is \(d<rk\), because \(\Gamma\) is
reduced. This contradiction proves \(\kappa(B)\geq1\).
\end{proof}

We next record an observation that will exclude rational fibrations on $S$.
\begin{lemma}\label{lem:vertical-branch}
Let \(f\colon S\to C\) be a fibration whose general fiber is
\(\prj{1}\), and let \(\pi\colon X\to S\) be a finite cover of
a smooth surface. If the branch divisor of \(\pi\) is contained in
finitely many fibers of \(f\), then \(X\) is birationally ruled.
\end{lemma}

\begin{proof}
For general \(c\in C\), the induced finite map
\(\pi^{-1}(f^{-1}(c))\to f^{-1}(c)\simeq\prj{1}\) is étale. Every finite
étale cover of \(\prj{1}\) is a disjoint union of copies of \(\prj{1}\).
In the Stein factorization of \(f\circ\pi\), the general connected
fiber is therefore \(\prj{1}\). Hence \(X\) is birationally ruled.
\end{proof}

\begin{lemma}\label{lem:isotropic-fibration}
In the cyclic-cover setup, suppose that \(P_B^2=0\).
Then \(P_B\) is semiample and determines a genus-one fibration \(f\colon S\to\prj{1}\). Every component of \(\Gamma\) is vertical, and \(K_S\cdot P_B=0\).
\end{lemma}

\begin{proof}
By \cref{lem:covering-Iitaka} and Zariski decomposition,
\(\kappa(P_B)=\kappa(B)\geq1\).
Take a divisible multiple \(nP_B\) whose linear system has nonzero
moving part \(M\). Its fixed part and \(M\) are orthogonal to
\(P_B\), by nefness and \(P_B^2=0\). The divisor \(M\) is nef and
\(M^2\geq0\). The Hodge index theorem therefore gives
\(M^2=0\) and \(M\equiv aP_B\), with \(a>0\).
Two general members of \(|M|\) have no common component and
intersection zero, so \(|M|\) is basepoint-free. On a rational
surface numerical and linear equivalence agree after clearing
denominators. Thus \(P_B\) is semiample and determines a fibration
with connected fibers; its base is \(\prj{1}\).

Since \(P_B\cdot\Gamma=rP_B^2+rP_B\cdot E=0\), the branch is
vertical. If \(Z\) is a general fiber, then
\[
 2g(Z)-2=K_S\cdot Z\leq0
\]
by anticanonical effectivity. Since genus zero would make \(X\) ruled by \cref{lem:vertical-branch}, we deduce that \(g(Z)=1\). Since \(Z^2=0\), adjunction gives \(K_S\cdot Z=0\), and \(Z\) is numerically a positive multiple of \(P_B\); hence \(K_S\cdot P_B=0\).
\end{proof}

\begin{lemma}[Positivity of the character adjoints]
\label{lem:invariant-nonnef}
In the cyclic-cover setup, assume that \(\rho(S)\geq3\).
Then for every $0\le j<r$, the divisor \(D_j\) is ample and basepoint-free, and
\begin{equation}\label{eq:Dj-curve}
 D_j\cdot C\geq p+2
 \qquad\text{for every irreducible curve }C\subset S.
\end{equation}
Moreover, for $F\in|-K_S|$
\[
 F\cdot D_j\geq p+2.
\]
Equality occurs only when \(F\) is integral, $F^2=1$,
and $B\cdot F>0$.
In particular, $F\cdot D\geq p+3$. Consequently, \(D\) satisfies \(\NP_p\).
\end{lemma}

\begin{proof}
Put \(Q=A-(p+3)F\) and \(C_j=Q+(r-1-j)B\).
The cone-theorem argument in \cref{prop:volume} shows that \(Q\) is
nef. If \(B\cdot C\geq0\), then \(C_j\cdot C\geq0\).
If \(B\cdot C=-n<0\), put \(g=g(C)\). By
\cref{prop:fractional-negative},
\[
 C_j\cdot C=A\cdot C-2(p+3)+2(p+3)g+\bigl((p+3)r-(r-1-j)\bigr)n
 \geq(r-1)(p+2)+j>0.
\]
Thus \(C_j\) is nef and
\begin{equation}\label{eq:Dj-decomp}
 D_j=C_j+(p+2)F.
\end{equation}
We check the required lower bound in three cases.

\smallskip\noindent
\emph{Case 1.} If $B\cdot C<0$, then
\[
 D_j\cdot C=A\cdot C+2g(C)-2+(j+1)n\geq p+2.
\]

\smallskip\noindent
\emph{Case 2.} If \(B\cdot C\geq0\) and \(F\cdot C\leq0\), then the bound follows
directly from \(D_j=A-F+(r-1-j)B\).

\smallskip\noindent
\emph{Case 3.} If \(B\cdot C\geq0\) and \(F\cdot C>0\), it follows
from \eqref{eq:Dj-decomp}. This proves \eqref{eq:Dj-curve}.

By \cref{prop:volume}, $D_j$ is ample and
basepoint-free. Since $F\in |-K_S|$ is nonzero and effective, $F\cdot D_j\ge p+2$.

We now determine when $F\cdot D_j$ can attain the value $p+2$.

If \(F\) has at least two irreducible components counted with
multiplicity, then
\(F\cdot D_j\geq2(p+2)\geq p+3\). Suppose that \(F\) is integral,
and put \(t=r-1-j\). If \(B\cdot F<0\), it is a smooth branch
component and \(F^2=rB\cdot F\). Hence
\[
 D_j\cdot F=A\cdot F+(t-r)B\cdot F\geq p+3.
\]
If \(B\cdot F\geq0\) and \(F^2\leq0\), the same bound is immediate.
If \(F^2\geq2\), use \eqref{eq:Dj-decomp} to obtain
\(D_j\cdot F\geq2(p+2)\).

It remains to consider \(F^2=1\) and \(B\cdot F\geq0\). Then
\[
 D_j\cdot F=A\cdot F-1+tB\cdot F\geq p+2.
\]
Equality implies \(A\cdot F=p+3\), so \(Q\cdot F=0\).
The Hodge index theorem and nefness of \(Q\) give \(Q\equiv0\),
hence \(A\equiv (p+3)F\). If \(B\cdot F=0\), then
\(A\cdot\Gamma=(p+3)rB\cdot F=0\), contradicting ampleness and
\(\Gamma\neq0\). Thus \(B\cdot F>0\), and equality can hold only
when \(t=0\). For \(D=D_0\), one has \(t=r-1>0\), so
\(D\cdot F\geq p+3\). Finally, the Gallego--Purnaprajna criterion
\cite[Theorem 1.3]{GP01} proves \(\NP_p\) for \(D\).
\end{proof}

\begin{proposition}\label{prop:rounded-adjunction}
In the cyclic-cover setup, for every integer \(0\leq t<r\),
\begin{equation*}
  K_S+A+tB\quad\text{is \(p\)-very ample}.
\end{equation*}
\end{proposition}

\begin{proof}
The case \(t=0\) follows from
\cref{lem:adjoint-separation}\textup{(i)}.
Assume henceforth that \(1\leq t<r\).

Keep notation in \cref{prop:fractional-negative}, and
set \(a_\alpha=A\cdot E_\alpha\) and
\[
 \lambda_{\alpha,t}
 =\max\left\{0,\frac{tn_\alpha+p+3-a_\alpha}{rn_\alpha}\right\},
 \qquad
 M_t=A+tB-\sum_\alpha\lambda_{\alpha,t}E_\alpha.
\]
Since \(a_\alpha\geq p+3\), one has
\(0\leq\lambda_{\alpha,t}\leq t/r<1\). Furthermore,
\[
 M_t=A+tP_B+
 \sum_\alpha\left(\frac tr-\lambda_{\alpha,t}\right)E_\alpha.
\]
If \(C\neq E_\alpha\) for every \(\alpha\), then
\(M_t\cdot C\geq A\cdot C\geq p+3\), while
\[
 M_t\cdot E_\alpha
 =a_\alpha-tn_\alpha+rn_\alpha\lambda_{\alpha,t}\geq p+3.
\]
Thus \(M_t\) is nef and has degree at least \(p+3\) on every curve.

Since \(rP_B\sim\Gamma-\sum E_\alpha\),
\[
 M_t-A\sim_{\mathbb{Q}}
 \frac tr\left(\Gamma-\sum_\alpha E_\alpha\right)
 +\sum_\alpha\left(\frac tr-\lambda_{\alpha,t}\right)E_\alpha
\]
is effective. Since \(A+M_t\) is ample,
\[
 M_t^2-A^2=(M_t-A)\cdot(M_t+A)\geq0.
\]
It follows that \(M_t^2\geq A^2\geq (p+3)^2\), and
\(\ceil{M_t}=A+tB\).

Let \(Z\) be a cluster of length \(p+1\). In
\cref{lem:Langer-rounded} take
\[
 \alpha=(p+2)^2,\qquad \beta=p+2.
\]
By \eqref{eq:Langer-cluster-bound}, \(\delta_Z\leq4(p+1)\), while
\[
 \alpha-4(p+1)=p^2\geq0,
 \qquad
 4\beta\left(1-\frac{\beta}{\alpha}\right)
 =4(p+2)\left(1-\frac1{p+2}\right)=4(p+1).
\]
Thus both inequalities in \eqref{eq:Langer-parameters} hold. Moreover,
\[
 M_t^2\geq(p+3)^2>\alpha,
 \qquad M_t\cdot C\geq p+3>\beta.
\]
Finally \(M_t\) was obtained from the integral divisor \(A+tB\) by
subtracting coefficients in \([0,1)\) along the pairwise disjoint smooth
curves \(E_\alpha\). Hence the fractional support of $M_t$ is simple normal
crossings and \(\ceil{M_t}=A+tB\). Therefore \cref{lem:Langer-rounded} implies that \(K_S+\ceil{M_t}=K_S+A+tB\) separates every cluster of length \(p+1\); that is, it is $p$-very ample.
\end{proof}

\section{Irregularity and regular covers}
\label{sec:regular}

For regular covers, the character decomposition gives the cohomology vanishing required by the weight-one criterion. We first show that a big positive part \(P_B\) forces \(X\) to be regular, and then obtain the specialty bound by restricting to an anticanonical divisor.

\begin{lemma}\label{lem:regularity}
In the cyclic-cover setup,
\begin{equation*}
 q(X)=\sum_{j=1}^{r-1}h^1(S,\sshf{S}(-jB)).
\end{equation*}
If \(P_B^2>0\), then \(q(X)=0\). In particular, every cover of
general type is regular.
\end{lemma}

\begin{proof}
The asserted equality follows from
\eqref{eq:cyclic-algebra} and \(H^1(S,\sshf{S})=0\).
If \(P_B\) is big, then \(jP_B\) is big and nef. Moreover,
\(\ceil{jP_B}=jB\) for \(1\leq j<r\), and its fractional support is a
disjoint union of smooth curves. Serre duality and Kawamata--Viehweg vanishing therefore give
\[
 h^1(S,\sshf{S}(-jB))=h^1(S,\sshf{S}(K_S+jB))=0.
\]
Thus \(q(X)=0\).

If \(X\) is of general type, the canonical bundle formula and
finiteness of \(\pi\) imply that \(K_S+(r-1)B\) is big. Adding the
effective divisor \(-K_S\) shows that \((r-1)B\), and hence \(B\),
is big. Hence \(P_B^2>0\), and $q(X)=0$.
\end{proof}

\begin{lemma}\label{lem:scalar-big}
Assume that \(\rho(S)\geq3\).
For every \(1\leq j<r\), the pair \((D,T=jB)\) satisfies
\eqref{eq:weight-hyp2}--\eqref{eq:weight-hyp3}.
\end{lemma}

\begin{proof}
By \cref{lem:invariant-nonnef}, $D$ is ample, so
\(D\cdot(jB)=(j/r)D\cdot\Gamma>0\).
Since \(rB\sim\Gamma\neq0\) and \(-K_S\) is nonzero effective,
we have
\[
 H^0(S,\sshf{S}(-jB))=H^0(S,\sshf{S}(K_S-jB))=0.
\]
By \cref{lem:invariant-nonnef} again, the divisor
\(qD-jB=D_j+(q-1)D\) is ample
for \(q\geq1\). Harbourne's theorem therefore gives
\(H^1(S,\sshf{S}(qD-jB))=0\).
\end{proof}

\begin{lemma}
\label{lem:anti-defect}
Assume that \(q(X)=0\) and \(\rho(S)\geq3\). For every
\(1\leq j<r\),
\begin{equation}\label{eq:anti-defect}
 (-K_S)\cdot D_j-h^1(S,\sshf{S}(jB))\geq p+2.
\end{equation}
\end{lemma}

\begin{proof}
Fix \(F\in|-K_S|\). By \cref{lem:anti-subcurves}, the curve \(F\)
is connected and Gorenstein, with \(\pa(F)=1\) and
\(\omega_F\simeq\sshf{F}\). The assumption that $q(X)=0$, \cref{lem:regularity} and Serre duality give
\(H^1(S,\sshf{S}(K_S+jB))=0\), while
\(H^0(S,\sshf{S}(-jB))=0\).
The exact sequence
\[\ses{\sshf{S}(K_S+jB)}{\sshf{S}(jB)}{\sshf{F}(jB)}\]
therefore yields
\begin{equation}\label{eq:h1-restrict-F}
 h^1(S,\sshf{S}(jB))=h^1(F,\sshf{F}(jB))
 =h^0(F,\sshf{F}(-jB)).
\end{equation}

There are two cases. Suppose first that \(F\cdot D_j\geq p+3\).
For every proper positive subcurve \(Y<F\),
\cref{lem:anti-subcurves,lem:invariant-nonnef} give
\(\pa(Y)\leq0\) and \(D_j\cdot Y\geq p+2\).
Consequently, for a cluster \(Z\subset F\) of length \(p+2\) and
every nonzero subcurve \(Y\subseteq F\),
\[
 D_j\cdot Y>2\pa(Y)-2+\operatorname{length}(Z\cap Y).
\]
The curve embedding theorem of Catanese--Franciosi--Hulek--Reid
\cite[Theorem 1.1]{CFHR} shows that \(D_j|_F\) separates \(Z\).
The same strict degree inequalities, automatic adjunction
\cite[Lemma 2.4]{CFHR}, and Serre duality give
\(H^1(F,\sshf{F}(D_j))=0\).
Riemann--Roch on \(F\) therefore gives
\begin{equation}\label{eq:h0-Dj-F}
 h^0(F,\sshf{F}(D_j))=D_j\cdot F.
\end{equation}

Since \(D\) is basepoint-free, a general section of \(D|_F\) is a
nonzerodivisor. Its zero scheme has length \(D\cdot F\geq p+3\);
choose within it a subcluster \(Z\) of length \(p+2\).
Multiplication by this section injects
\[
 H^0(F,\sshf{F}(-jB))\hookrightarrow
 H^0(F,\is{Z}\otimes\sshf{F}(D_j)).
\]
The target has codimension \(p+2\) in
\(H^0(F,\sshf{F}(D_j))\), by separation
of \(Z\). Combining \eqref{eq:h1-restrict-F} and
\eqref{eq:h0-Dj-F} proves \eqref{eq:anti-defect}.

The remaining case is \(D_j\cdot F=p+2\).
Then by \cref{lem:invariant-nonnef}, the curve \(F\) is integral with
\(B\cdot F>0\). Thus \(-jB|_F\) has negative degree and no nonzero
section. Equation \eqref{eq:h1-restrict-F} gives
\(h^1(S,\sshf{S}(jB))=0\),
which proves the desired inequality in this case too.
\end{proof}

\begin{proposition}\label{prop:regular-case}
In the cyclic-cover setup, if \(q(X)=0\), then \(L\) satisfies
\(\NP_p\).
\end{proposition}

\begin{proof}
The case \(\rho(S)\leq2\) is proved in
\cref{prop:low-rank}. Assume therefore that \(\rho(S)\geq3\).
By \cref{lem:invariant-nonnef}, the divisor \(D\) satisfies \(\NP_p\).
For \(1\leq j<r\), set \(T=jB\).
Since \(D_j=K_S+A+(r-1-j)B\),
\cref{prop:rounded-adjunction} shows that \(D_j\)
is \(p\)-very ample. \Cref{lem:scalar-big,lem:anti-defect} verify the remaining hypotheses of \cref{cor:cyclic-reduction}, which gives the assertion.
\end{proof}

By \cref{lem:regularity,prop:regular-case}, \(L\) satisfies
\(\NP_p\) whenever \(P_B^2>0\); in particular, this holds whenever
\(X\) is of general type.

The hypothesis \(P_B^2>0\) in \cref{lem:regularity} cannot be dropped: when \(P_B^2=0\), the cover can have arbitrarily large irregularity.

\begin{example}\label[example]{ex:large-irregularity}
Let \(f\colon S\to\prj{1}\) be a rational elliptic surface with fiber
\(F=-K_S\) (see \cref{sec:elliptic}); for instance, the blow-up of
\(\prj{2}\) at the nine base points of a general pencil of cubics.
For \(g\geq1\), let \(\Gamma\) be the union of \(2g+2\) distinct smooth
fibers, and let \(\pi\colon X\to S\) be the double cover branched
along \(\Gamma\), so that \(B=(g+1)F\). Then \(X=S\times_{\prj{1}}C\),
where \(C\) is the double cover of $\prj{1}$ branched over
\(f(\Gamma)\), hence is of genus \(g\). Here \(P_B^2=B^2=0\), and
\(K_X=\pi^*(gF)\) is effective, so \(X\) is nonruled. By \cref{lem:regularity} and Serre duality,
\[
 q(X)=h^1(S,\sshf{S}(-(g+1)F))=g, \qquad
 h^1(S,\sshf{S}(B))=h^1(S,\sshf{S}(-(g+2)F))=g+1.
\]
Since \(D_1=K_S+A\), one has \((-K_S)\cdot D_1=F\cdot A\). So for fixed
\(A\) and \(g\gg0\), the inequality \eqref{eq:anti-defect} fails for
\(j=1\), and \cref{prop:weight-one} does not apply to \((D,B)\). Such
covers are treated in \cref{sec:elliptic} by a different method.
\end{example}

\section{Rational elliptic surfaces}\label{sec:elliptic}

Recall that a \emph{Halphen surface of index \(\iota\geq1\)} is a smooth rational surface \(S\) such that \(|-\iota K_S|\) is a basepoint-free pencil. Its members are the fibers of a relatively minimal genus-one fibration \(S\to\prj{1}\), and every relatively minimal genus-one fibration on a smooth rational surface arises in this way; see \cite[Section~2.2]{CantatDolgachev}. The Halphen surfaces of index one are the rational elliptic surfaces with a section.

Throughout this section, we assume that \(|-K_S|\) is a basepoint-free pencil, so that \(S\) is a Halphen surface of index one. We write \(f\colon S\to\prj{1}\) for the elliptic fibration and \(F=-K_S\) for the fiber class; thus \(F^2=0\) and \(\rho(S)=10\). Two features of this situation drive the arguments: a general fiber is a smooth elliptic curve, on which kernel bundles are semistable, so that the required vanishing on fibers follows from slope considerations (\cref{lem:fiber-restriction}\textup{(ii)}); and the fibers move, so that \(F\) can be added to a divisor one copy at a time. In \cref{sec:add-fibers} we use this to treat nef covering data. Non-nef covering data are treated in \cref{sec:relative} with relative kernel bundles, which change only by a twist when fibers are added.

The following lemma collects what we need about restriction to fibers.

\begin{lemma}\label{lem:fiber-restriction}    
\begin{enumerate} 
\item Let \(H\) be a globally generated divisor on \(S\) with \(H^1(S,\sshf{S}(H-F))=0\), and let \(F_0\) be a fiber of \(f\). Then the restriction map \(H^0(S,\sshf{S}(H))\to H^0(F_0,\sshf{F_0}(H))\) is surjective, and there is an exact sequence
\begin{equation}\label{eq:complete-restriction-kernel}
 0\longrightarrow H^0(S,\sshf{S}(H-F))\otimes\sshf{F_0}
 \longrightarrow M_H|_{F_0}\longrightarrow M_{H|_{F_0}}\longrightarrow0.
\end{equation}
\item Let \(C\) be a smooth elliptic curve, let \(J_1,\dots,J_k\) be line bundles of degree \(e\geq2\) on \(C\), let \(w_1,\dots,w_k\geq0\), and let \(T\) be a line bundle on \(C\) with \(\deg T>(w_1+\cdots+w_k)\,e/(e-1)\). Then
\[
 H^1\Bigl(C,\textstyle\bigotimes_{i=1}^k\bigwedge^{w_i}M_{J_i}\otimes T\Bigr)=0.
\]
\end{enumerate}
\end{lemma}

\begin{proof}
(i) Since \(F_0\sim F\), the restriction map is surjective, with kernel \(H^0(S,\sshf{S}(H-F))\) embedded by multiplication with an equation of \(F_0\). In particular \(\sshf{F_0}(H)\) is globally generated. Restrict the evaluation sequence of \(H\) to \(F_0\), and compare it with the evaluation sequence of \(\sshf{F_0}(H)\) via the restriction map; the snake lemma gives \eqref{eq:complete-restriction-kernel}.

(ii) By \cite{Butler}, each \(M_{J_i}\) is stable of rank \(e-1\) and degree \(-e\). In characteristic zero, exterior powers and tensor products of semistable bundles are semistable. Hence the bundle in question is zero or semistable of slope \(\deg T-(w_1+\cdots+w_k)e/(e-1)>0\). By Serre duality, such a bundle has no \(H^1\).
\end{proof}

\subsection{Adding fibers and nef covering data}\label{sec:add-fibers}

\begin{lemma}\label{lem:add-anticanonical}
Let \(H\) be an ample divisor on \(S\) such that \(H-F\) is nef, and put \(e=H\cdot F\). If \(e\geq p+3\), then
\begin{equation}\label{eq:add-anticanonical}
 H^1\!\left(S,\bigwedge^hM_{H+cF}\otimes\sshf{S}(\ell H+dF)\right)=0
\end{equation}
for \(c\geq0\), \(\ell\geq1\), \(0\leq h\leq p+1\), and \(d\geq h\).
\end{lemma}

\begin{proof}
Put \(H_t=H+tF\). Since \(H-F\) is nef and \((H-F)\cdot F=e\geq2\), Harbourne's theorem shows that \(H-F\) and every \(H_t\) are globally generated with vanishing \(H^1\), and that \(H^1(S,\sshf{S}(\ell H+tF))=0\) for \(\ell\geq1\) and \(t\geq0\). In particular, \cref{lem:fiber-restriction}\textup{(i)} applies to every \(H_t\). Fix a smooth fiber, which we also denote by \(F\), and let \(i\colon F\hookrightarrow S\) be the inclusion. Note that \(\sshf{F}(H_t)\) has degree \(e\) for every \(t\), since \(F^2=0\), and that \(p+1\leq e-2\).

\smallskip\noindent\emph{Step 1: the case \(c=0\).} We show that
\begin{equation}\label{eq:initial-kernel-twists}
 H^1\Bigl(S,\textstyle\bigwedge^aM_H\otimes\sshf{S}(\ell H+tF)\Bigr)=0
 \qquad(0\leq a\leq p+1,\ \ell\geq1,\ t\geq0).
\end{equation}
For \(t=0\), this is the kernel-bundle form of Property \(\NP_p\) for \(H\), which holds by \cite[Theorem 1.3]{GP01} because \(H\cdot(-K_S)=e\geq p+3\). Let \(t\geq1\). By \eqref{eq:complete-restriction-kernel}, the restriction of \(\bigwedge^aM_H\) to \(F\) has a filtration whose graded pieces are the bundles \(\bigwedge^{a'}M_{H|_F}\), \(a'\leq a\), tensored with trivial bundles. The line bundle \(\sshf{F}(\ell H+tF)\) has degree \(\ell e>ae/(e-1)\), so \cref{lem:fiber-restriction}\textup{(ii)} shows that the restriction to \(F\) of the bundle in \eqref{eq:initial-kernel-twists} has no \(H^1\). The restriction sequence therefore passes from \(t-1\) to \(t\).

\smallskip\noindent\emph{Step 2: comparing \(H_{t-1}\) and \(H_t\).} Let \(t\geq1\) and put \(W_t=H^0(F,\sshf{F}(H_t))\). By \cref{lem:fiber-restriction}\textup{(i)}, the sequence
\[
 0\longrightarrow H^0(S,\sshf{S}(H_{t-1}))\longrightarrow H^0(S,\sshf{S}(H_t))\longrightarrow W_t\longrightarrow0
\]
is exact. Tensor it with \(\sshf{S}\), and map it by evaluation to the exact sequence \(0\to\sshf{S}(H_{t-1})\to\sshf{S}(H_t)\to i_*\sshf{F}(H_t)\to0\). The first vertical map is surjective, since \(H_{t-1}\) is globally generated, so the snake lemma gives an exact sequence
\begin{equation}\label{eq:elementary-kernel-comparison}
 0\longrightarrow M_{H_{t-1}}\longrightarrow M_{H_t}
 \longrightarrow\shf{V}_t\longrightarrow0,
 \qquad
 \shf{V}_t=\ker\bigl(W_t\otimes\sshf{S}\longrightarrow i_*\sshf{F}(H_t)\bigr).
\end{equation}
The sheaf \(\shf{V}_t\) is locally free of rank \(e\). Indeed, near a point of \(F\), let \(x\) be a local equation of \(F\), and choose a local frame \(u_1,\dots,u_e\) of \(W_t\otimes\sshf{S}\) such that \(u_1\) maps to a generator of \(\sshf{F}(H_t)\) and \(u_2,\dots,u_e\) restrict to a frame of \(M_{H_t|_F}\). Then \(\shf{V}_t\) is generated by \(xu_1,u_2,\dots,u_e\). For \(v\geq1\), the same frame shows that
\begin{equation}\label{eq:elementary-exterior}
 0\longrightarrow\bigwedge^vW_t\otimes\sshf{S}(-F)
 \longrightarrow\bigwedge^v\shf{V}_t
 \longrightarrow i_*\bigwedge^vM_{H_t|_F}\longrightarrow0
\end{equation}
is exact: the image of the first map is \(x\bigwedge^v(W_t\otimes\sshf{S})\), and the quotient is generated by the wedge products of \(u_2,\dots,u_e\). Similarly, restricting \(\shf{V}_t\) to \(F\) gives
\begin{equation}\label{eq:elementary-on-F}
 0\longrightarrow\sshf{F}(H_t-F)\longrightarrow\shf{V}_t|_F\longrightarrow M_{H_t|_F}\longrightarrow0.
\end{equation}

\smallskip\noindent\emph{Step 3: filtering.} Iterating \eqref{eq:elementary-kernel-comparison} for \(t=1,\dots,c\) filters \(M_{H+cF}=M_{H_c}\) with graded pieces \(M_H,\shf{V}_1,\dots,\shf{V}_c\), and hence filters \(\bigwedge^hM_{H+cF}\) with graded pieces
\[
 \bigwedge^aM_H\otimes\bigwedge^{v_1}\shf{V}_1\otimes\cdots\otimes\bigwedge^{v_c}\shf{V}_c,
 \qquad a+v_1+\cdots+v_c=h.
\]
Twist such a piece by \(\ell H+dF\), and apply \eqref{eq:elementary-exterior} to the factors with \(v_t\geq1\), one at a time. Since the remaining factors are locally free, every sequence obtained in this way is exact. At each step we either replace \(\bigwedge^{v_t}\shf{V}_t\) by \(\bigwedge^{v_t}W_t\otimes\sshf{S}(-F)\), or pass to the quotient term, which is a sheaf on \(F\); in the latter case all remaining factors are restricted to \(F\). The resulting graded pieces are of two kinds.
\begin{enumerate}
\item[(a)] Pieces on \(S\). Up to a vector-space factor, these are \(\bigwedge^aM_H\otimes\sshf{S}(\ell H+d'F)\) with \(d'\geq d-(h-a)\geq0\), since at most \(h-a\) copies of \(F\) are subtracted and \(d\geq h\). They have no \(H^1\), by \eqref{eq:initial-kernel-twists}.
\item[(b)] Pieces on \(F\). By \eqref{eq:complete-restriction-kernel} and \eqref{eq:elementary-on-F}, these are filtered by bundles \(\bigwedge^{w_0}M_{H|_F}\otimes\bigotimes_t\bigwedge^{w_t}M_{H_t|_F}\otimes T\), where \(w_0+\sum_tw_t\leq h\) and \(T\) is a line bundle of degree at least \(\ell e\). As \(\ell e>he/(e-1)\), they have no \(H^1\), by \cref{lem:fiber-restriction}\textup{(ii)}.
\end{enumerate}
The long exact cohomology sequences of these filtrations prove \eqref{eq:add-anticanonical}.
\end{proof}

\begin{proposition}[Nef covering data]\label{prop:elliptic-nef}
In the cyclic-cover setup, suppose that \(|-K_S|\) is a basepoint-free pencil and that \(B\) is nef. Then \(L\) satisfies \(\NP_p\).
\end{proposition}

\begin{proof}
If \(B^2>0\), then \(P_B=B\) is big, and the assertion follows from \cref{lem:regularity,prop:regular-case}. Assume that \(B^2=0\). Then \(P_B=B\), and \cref{lem:isotropic-fibration} gives \(B\cdot F=0\). Since \(B^2=F^2=0\) and \(B\not\equiv0\), the Hodge index theorem gives \(B\equiv bF\) with \(b>0\). A \((-1)\)-curve \(E_0\) on \(S\) satisfies \(F\cdot E_0=1\) by adjunction, so \(b=B\cdot E_0\) is an integer. Numerical and linear equivalence agree on \(S\), so \(B\sim bF\).

Put
\[
 Q=A-(p+3)F,\qquad H=Q+F=A-(p+2)F,\qquad
 \alpha=(r-1)b,\qquad c=p+\alpha+1.
\]
The cone-theorem argument in \cref{prop:volume} shows that \(Q\) is nef. The divisor \(H\) is positive on every curve \(C\): if \(F\cdot C=0\), then \(H\cdot C=A\cdot C\geq p+3\), and otherwise \(H\cdot C\geq F\cdot C>0\). Since \(H^2=Q^2+2Q\cdot F>0\), the divisor \(H\) is ample. Moreover, \(H-F=Q\) is nef and \(H\cdot F=A\cdot F\geq p+3\), so \cref{lem:add-anticanonical} applies to \(H\). As \(K_S=-F\) and \(B\sim bF\),
\[
 D=H+cF,\qquad D_j=H+(c-jb)F.
\]
For \(1\leq j<r\), \(q\geq0\), and \(0\leq i\leq p\), put \(\ell=q+1\) and \(h=i+1\). Then
\[
 D_j+qD=\ell H+dF,\qquad
 d=\ell c-jb\geq c-\alpha=p+1\geq h,
\]
so \eqref{eq:add-anticanonical} gives \(H^1\bigl(S,\bigwedge^{i+1}M_D\otimes\sshf{S}(D_j+qD)\bigr)=0\). By \cref{lem:kernel-criterion}, condition \textup{(ii)} of \cref{thm:transfer} holds for every nontrivial character. Finally, \(D\) is ample and globally generated with \(D\cdot F=A\cdot F\geq p+3\), so it satisfies \(\NP_p\) by \cite[Theorem 1.3]{GP01}; this is condition \textup{(i)}, and condition \textup{(iii)} is \cref{prop:high-weights}. \Cref{thm:transfer} gives \(\NP_p\) for \(L\).
\end{proof}

\subsection{Relative kernel bundles and non-nef covering data}\label{sec:relative}

For a divisor \(H\) on \(S\), the \emph{relative kernel bundle} of \(H\) is
\[
 M_{H/f}=\ker\bigl(f^*f_*\sshf{S}(H)\longrightarrow\sshf{S}(H)\bigr),
\]
the kernel of the relative evaluation map. Since \(f^*\sshf{\prj{1}}(1)=\sshf{S}(F)\), the projection formula gives \(f_*\sshf{S}(H+mF)=f_*\sshf{S}(H)\otimes\sshf{\prj{1}}(m)\), and hence
\begin{equation}\label{eq:relative-twist}
 M_{(H+mF)/f}\simeq M_{H/f}\otimes\sshf{S}(mF)\qquad(m\in\mathbb Z).
\end{equation}
This is why we pass to relative kernel bundles: adding fibers to \(H\) changes \(M_{H/f}\) only by a twist, whereas it changes \(M_H\) as in \eqref{eq:elementary-kernel-comparison}.

\begin{lemma}[From complete to relative kernel bundles]\label{lem:global-relative}
Let \(H\) be an ample, globally generated divisor on \(S\) with \(H^1(S,\sshf{S}(H-F))=0\) and \(e=H\cdot F\geq2\).
\begin{enumerate}
\item The sheaf \(\shf{E}=f_*\sshf{S}(H)\) is a globally generated vector bundle on \(\prj{1}\), and there is an exact sequence
\begin{equation}\label{eq:complete-relative}
 0\longrightarrow\sshf{S}(-F)^{\oplus N}\longrightarrow M_H\longrightarrow M_{H/f}\longrightarrow0,
 \qquad N=\deg\shf{E}.
\end{equation}
\item Let \(v\geq1\) with \(e\geq v+2\), and let \(G\) be an ample divisor with \(G\cdot F=e\). If
\begin{equation}\label{eq:mixed-initial}
 H^1\Bigl(S,\textstyle\bigwedge^wM_H\otimes\sshf{S}(G+qH)\Bigr)=0
 \qquad(0\leq w\leq v,\ q\geq0),
\end{equation}
then
\begin{equation}\label{eq:mixed-relative}
 H^1\Bigl(S,\textstyle\bigwedge^vM_{H/f}\otimes\sshf{S}(G+qH+kF)\Bigr)=0
 \qquad(q\geq0,\ k\geq1).
\end{equation}
\end{enumerate}
\end{lemma}

\begin{proof}
(i) By Harbourne's theorem, \(H^1(S,\sshf{S}(H))=0\); and \(H^2(S,\sshf{S}(H-F))
=0\) by Serre duality. Restricting to a fiber \(F_t\sim F\), we obtain \(H^1(F_t,\sshf{F_t}(H))=0\) for every \(t\in\prj{1}\); and by \cref{lem:fiber-restriction}\textup{(i)}, the restriction maps \(H^0(S,\sshf{S}(H))\to H^0(F_t,\sshf{F_t}(H))\) are surjective. By cohomology and base change, \(\shf{E}\) is locally free with fibers \(H^0(F_t,\sshf{F_t}(H))\), and it is globally generated, since \(H^0(\prj{1},\shf{E})=H^0(S,\sshf{S}(H))\). The evaluation map of \(H\) factors as
\[
 H^0(S,\sshf{S}(H))\otimes\sshf{S}\longrightarrow f^*\shf{E}\longrightarrow\sshf{S}(H).
\]
The first map is the pullback of the evaluation map of \(\shf{E}\); it is surjective, with kernel \(f^*M_{\shf{E}}\). The composite is surjective because \(H\) is globally generated, hence so is the second map, the relative evaluation map, whose kernel is \(M_{H/f}\). This gives an exact sequence \(0\to f^*M_{\shf{E}}\to M_H\to M_{H/f}\to0\). Finally, write \(\shf{E}\simeq\bigoplus_\nu\sshf{\prj{1}}(d_\nu)\); global generation gives \(d_\nu\geq0\). Since the kernel bundle of \(\sshf{\prj{1}}(d)\) is \(\sshf{\prj{1}}(-1)^{\oplus d}\), we get \(f^*M_{\shf{E}}\simeq\sshf{S}(-F)^{\oplus N}\) with \(N=\sum_\nu d_\nu=\deg\shf{E}\).

(ii) We first show that
\begin{equation}\label{eq:mixed-initial-t}
 H^1\Bigl(S,\textstyle\bigwedge^wM_H\otimes\sshf{S}(G+qH+tF)\Bigr)=0
 \qquad(0\leq w\leq v,\ q,t\geq0).
\end{equation}
For \(t=0\) this is \eqref{eq:mixed-initial}. To pass from \(t-1\) to \(t\), restrict to a smooth fiber \(F_0\). By \eqref{eq:complete-restriction-kernel}, \(\bigwedge^wM_H|_{F_0}\) is filtered by the bundles \(\bigwedge^{w'}M_{H|_{F_0}}\), \(w'\leq w\), tensored with trivial bundles, and \(\sshf{F_0}(G+qH+tF)\) has degree \((q+1)e>we/(e-1)\), since \(w\leq v\leq e-2\). Hence \cref{lem:fiber-restriction}\textup{(ii)} applies.

Next we show that
\begin{equation}\label{eq:mixed-H2}
 H^2\Bigl(S,\textstyle\bigwedge^{v-1}M_H\otimes\sshf{S}(G+qH+tF)\Bigr)=0
 \qquad(q,t\geq0).
\end{equation}
Put \(\Xi=G+qH+tF\). As \(\Xi\) is ample, \(H^2(S,\sshf{S}(\Xi))\simeq H^0(S,\sshf{S}(K_S-\Xi))^\vee=0\), which is \eqref{eq:mixed-H2} for \(v=1\). For \(v\geq2\), twist \eqref{eq:koszul-wedge}, with \(D=H\) and \(w=v-1\), by \(\Xi\). In the resulting cohomology sequence, \(H^2\bigl(S,\bigwedge^{v-1}M_H\otimes\sshf{S}(\Xi)\bigr)\) lies between \(H^1\bigl(S,\bigwedge^{v-2}M_H\otimes\sshf{S}(\Xi+H)\bigr)\), which vanishes by \eqref{eq:mixed-initial-t}, and \(\bigwedge^{v-1}H^0(S,\sshf{S}(H))\otimes H^2(S,\sshf{S}(\Xi))=0\).

Finally, let \(\shf{K}_v\) be the kernel of the surjection \(\bigwedge^vM_H\to\bigwedge^vM_{H/f}\) induced by \eqref{eq:complete-relative}. It is the image of the wedge-product map \(\sshf{S}(-F)^{\oplus N}\otimes\bigwedge^{v-1}M_H\to\bigwedge^vM_H\). Twist by \(G+qH+kF\), with \(k\geq1\). The source then has no \(H^2\), by \eqref{eq:mixed-H2} with \(t=k-1\); since \(S\) is a surface, neither has \(\shf{K}_v\otimes\sshf{S}(G+qH+kF)\). The cohomology sequence of \(0\to\shf{K}_v\to\bigwedge^vM_H\to\bigwedge^vM_{H/f}\to0\), twisted by \(G+qH+kF\), together with \eqref{eq:mixed-initial-t} for \(w=v\), gives \eqref{eq:mixed-relative}.
\end{proof}

\begin{lemma}[Non-nef covering data]\label{lem:small-vertical}
In the cyclic-cover setup, suppose that \(|-K_S|\) is a basepoint-free pencil, that \(P_B^2=0\), and that \(B\) is not nef. Then \(r=2\), and
\[
 B\sim_{\mathbb{Q}}bF+\tfrac12R,\qquad R=E_1+\cdots+E_s,
\]
where \(E_1,\dots,E_s\) are the curves on which \(B\) is negative; they are pairwise disjoint \((-2)\)-curves with \(F\cdot E_\alpha=0\), \(s\in\{4,8\}\), and \(2b\) is a positive integer. Put \(n=\floor b\), \(b_0=b-n\in\{0,\tfrac12\}\), and \(B_0=B-nF\). Then
\[
 H^0(S,\sshf{S}(-B_0))=H^0(S,\sshf{S}(K_S-B_0))=H^0(S,\sshf{S}(K_S+B_0))=0,
 \qquad h^1(S,\sshf{S}(B_0))\leq1.
\]
\end{lemma}

\begin{proof}
Let \(E_\alpha\) be a curve with \(B\cdot E_\alpha=-n_\alpha<0\). By \cref{prop:fractional-negative} and adjunction, \(F\cdot E_\alpha=2-2g(E_\alpha)-rn_\alpha\). Since \(F\) is nef, this forces \(g(E_\alpha)=0\), \(r=2\), \(n_\alpha=1\), and \(F\cdot E_\alpha=0\); thus \(E_\alpha\) is a \((-2)\)-curve, and \(B=P_B+\tfrac12R\) by \eqref{eq:zariski-B}. By \cref{lem:isotropic-fibration}, \(P_B\cdot F=0\). Since \(P_B^2=F^2=0\) and \(P_B\not\equiv0\) (\cref{prop:fractional-negative}), the Hodge index theorem gives \(P_B\equiv bF\) with \(b>0\). If \(E\) is a \((-1)\)-curve on \(S\), then \(F\cdot E=1\) by adjunction, so \(2b=(2B-R)\cdot E\) is an integer. As numerical and \(\mathbb Q\)-linear equivalence agree on \(S\), we get \(B\sim_{\mathbb Q}bF+\tfrac12R\) and \(2B_0\sim2b_0F+R\).

Since \(F^2=F\cdot R=0\) and \(R^2=-2s\), the formula for double covers gives
\[
 \chi(\sshf{X})=2\chi(\sshf{S})+\tfrac12B\cdot(B+K_S)=2+\tfrac18R^2=2-\tfrac{s}{4}.
\]
As \(X\) is not ruled, \(\chi(\sshf{X})\geq0\) by the classification of surfaces. Since \(s>0\), we conclude that \(s\in\{4,8\}\).

Both \(B_0\equiv b_0F+\tfrac12R\) and \(F+B_0\) are numerically nonzero effective \(\mathbb Q\)-divisors, so \(-B_0\) and \(K_S-B_0=-(F+B_0)\) are not effective. If an effective divisor \(V\) were linearly equivalent to \(K_S+B_0\), then \(V\cdot E_\alpha=-1\) would force \(V\geq R\); but \(V-R\equiv(b_0-1)F-\tfrac12R\) has negative degree on an ample divisor. Hence \(H^0(S,\sshf{S}(K_S+B_0))=0\).

Suppose that \(s=4\). Riemann--Roch gives \(\chi(\sshf{S}(-B_0))=1+\tfrac18R^2=0\), and \(h^0(S,\sshf{S}(-B_0))=h^2(S,\sshf{S}(-B_0))=0\) by the vanishings just proved; hence \(H^1(S,\sshf{S}(-B_0))=0\). Let \(F_0\) be a smooth fiber. In the cohomology sequence of
\[
 0\longrightarrow\sshf{S}(K_S+B_0)\longrightarrow\sshf{S}(B_0)\longrightarrow\sshf{F_0}(B_0)\longrightarrow0,
\]
the groups 
\[H^1(S,\sshf{S}(K_S+B_0))\simeq H^1(S,\sshf{S}(-B_0))^\vee,  H^2(S,\sshf{S}(K_S+B_0))\simeq H^0(S,\sshf{S}(-B_0))^\vee\] vanish. Hence \(h^1(S,\sshf{S}(B_0))=h^1(F_0,\sshf{F_0}(B_0))\leq1\), since \(\sshf{F_0}(B_0)\) has degree zero on the elliptic curve \(F_0\).

Suppose that \(s=8\). We claim that \(B_0\) is not effective. If \(V\sim B_0\) were effective, then \(V\cdot E_\alpha=-1\) would again force \(V\geq R\). If \(b_0=0\), then \(V-R\equiv-\tfrac12R\), which is impossible. If \(b_0=\tfrac12\), then \(2(V-R)+R\sim2B_0-R\sim F\), so \(2(V-R)+R\) is a fiber of \(f\) containing the eight pairwise disjoint curves \(E_\alpha\). But a fiber of a rational elliptic surface has at most nine components, and none of the fibers with nine components, of types \(\mathrm I_9\), \(\mathrm I_4^*\) and \(\mathrm{II}^*\), contains eight pairwise disjoint ones; see \cite[Table~(1.1)]{Miranda}. Thus \(H^0(S,\sshf{S}(B_0))=0\). Since \(\chi(\sshf{S}(B_0))=1+\tfrac18R^2=-1\) and \(H^2(S,\sshf{S}(B_0))\simeq H^0(S,\sshf{S}(K_S-B_0))^\vee=0\), we get \(h^1(S,\sshf{S}(B_0))=1\).
\end{proof}

\begin{theorem}\label{thm:elliptic-complete}
In the cyclic-cover setup, if \(|-K_S|\) is a basepoint-free pencil, then \(L\) satisfies \(\NP_p\).
\end{theorem}

\begin{proof}
By \cref{lem:regularity,prop:regular-case} we may assume that \(P_B^2=0\), and by \cref{prop:elliptic-nef} that \(B\) is not nef. We use the notation of \cref{lem:small-vertical}; in particular \(r=2\). Put
\[
 Q=A-(p+3)F,\qquad G=K_S+A,\qquad D^0=G+B_0,\qquad e=A\cdot F\geq p+3.
\]
Then \(D=G+B=D^0+nF\) and \(D_1=D-B=G\), so the unique nontrivial character module is \(N_{D,B}=\bigoplus_{q\geq0}H^0(S,\sshf{S}(G+qD))\). Since \(\rho(S)=10\), \cref{lem:invariant-nonnef} shows that \(D\) satisfies \(\NP_p\), and \cref{prop:high-weights} gives condition \textup{(iii)} of \cref{thm:transfer}. By \cref{lem:kernel-criterion}, it remains to prove that
\begin{equation}\label{eq:target-vanishing}
 H^1\Bigl(S,\textstyle\bigwedge^{i+1}M_D\otimes\sshf{S}(G+qD)\Bigr)=0
 \qquad(0\leq i\leq p,\ q\geq0).
\end{equation}

Since \(B\) is not a multiple of \(F\), \cref{lem:add-anticanonical} does not apply to \(D\). Instead, Step~1 proves the analogous vanishing for the divisors \(D^0-uF\) by means of the weight-one criterion and Castelnuovo--Mumford regularity; this settles the case \(n=0\). When \(n\geq1\), Step~3 compares \(M_D\) with its relative kernel bundle, which by \eqref{eq:relative-twist} is a twist of the relative kernel bundle of \(D^0-uF\).

We shall use that \(Q\) and \(Q+B_0\) are nef. For \(Q\), this is the cone-theorem argument in \cref{prop:volume}. For \(Q+B_0\), note that \((Q+B_0)\cdot E_\alpha=A\cdot E_\alpha-1\geq p+2\), and that \((Q+B_0)\cdot C\geq Q\cdot C\geq0\) for every other irreducible curve \(C\), because \(B_0\equiv b_0F+\tfrac12R\).

\smallskip\noindent\emph{Step 1: the divisors \(D^0-uF\).} Fix \(1\leq v\leq p+1\), put \(u=p+1-v\), and let
\[
 A_u=A-uF,\qquad H=D^0-uF,\qquad G_u=G-uF=K_S+A_u.
\]
Then \(A_u=Q+(v+2)F\) has degree at least \(v+2\) on every curve, and
\[
 H=Q+B_0+(v+1)F,\qquad G_u=Q+(v+1)F.
\]
Hence \(H\), \(H-F\), and \(G_u\) are positive on every curve and have positive square, so they are ample; by Harbourne's theorem they are globally generated with vanishing \(H^1\). Moreover, \(H\cdot F=G_u\cdot F=e\geq3\), so \(H\) is very ample by \cite[Theorem 1.3]{GP01}. We claim that
\begin{equation}\label{eq:small-data-kernels}
 H^1\Bigl(S,\textstyle\bigwedge^wM_H\otimes\sshf{S}(G_u+qH)\Bigr)=0
 \qquad(0\leq w\leq v,\ q\geq0).
\end{equation}
This follows from \cref{lem:weight-one-regularity}, applied with \((D,T,m)=(H,B_0,v-1)\), so that \(D-T=G_u\). Its hypotheses hold: \(G_u=K_S+A_u\) is \((v-1)\)-very ample by \cref{lem:adjoint-separation}\textup{(i)}; \(H\cdot B_0>0\); \(H^0(S,\sshf{S}(-B_0))=H^0(S,\sshf{S}(K_S-B_0))=H^0(S,\sshf{S}(K_S+B_0))=0\) by \cref{lem:small-vertical}; \(H^1(S,\sshf{S}(qH-B_0))=H^1(S,\sshf{S}(G_u+(q-1)H))=0\) for \(q\geq1\), by Harbourne's theorem; and
\[
 (-K_S)\cdot G_u-h^1(S,\sshf{S}(B_0))\geq e-1\geq p+2\geq v+1.
\]

\smallskip\noindent\emph{Step 2: the case \(n=0\).} Then \(D=D^0\). Taking \(v=p+1\), hence \(u=0\), in Step~1 gives \(H=D\) and \(G_u=G\), and \eqref{eq:small-data-kernels} becomes \eqref{eq:target-vanishing}.

\smallskip\noindent\emph{Step 3: the case \(n\geq1\).} The divisor \(D=D^0+nF\) is ample and globally generated, and \(D-F=D^0+(n-1)F\) is ample with \((D-F)\cdot F=e\geq3\), so \(H^1(S,\sshf{S}(D-F))=0\) by Harbourne's theorem. By \cref{lem:global-relative}\textup{(i)}, there is an exact sequence \(0\to\sshf{S}(-F)^{\oplus N}\to M_D\to M_{D/f}\to0\). Consequently, for \(i\leq p\) and \(q\geq0\), the bundle \(\bigwedge^{i+1}M_D\otimes\sshf{S}(G+qD)\) has a filtration whose graded pieces are direct sums of the bundles
\begin{equation}\label{eq:elliptic-mixed-filtration}
 \textstyle\bigwedge^vM_{D/f}\otimes\sshf{S}(G+qD-aF),\qquad a+v=i+1\leq p+1.
\end{equation}
It suffices to show that each of them has vanishing \(H^1\).

If \(v=0\), then \(a=i+1\leq p+1\), and \(G+qD-aF=Q+(p+2-a)F+qD\) is ample with anticanonical degree at least \(e\), so Harbourne's theorem applies.

If \(v\geq1\), let \(u=p+1-v\) and \(H=D^0-uF\) be as in Step~1. Since \(D=H+(n+u)F\), \eqref{eq:relative-twist} gives \(M_{D/f}\simeq M_{H/f}((n+u)F)\), and the bundle \eqref{eq:elliptic-mixed-filtration} becomes
\[
 \textstyle\bigwedge^vM_{H/f}\otimes\sshf{S}(G_u+qH+kF),\qquad
 k=(q+v)n+(q+v+1)u-a.
\]
As \(a=i+1-v\leq u\) and \(n\geq1\), we have \(k\geq(q+v)(n+u)\geq1\). By \eqref{eq:small-data-kernels}, \cref{lem:global-relative}\textup{(ii)} applies to \(H\) and \(G_u\), and shows that this bundle has vanishing \(H^1\). This proves \eqref{eq:target-vanishing}.
\end{proof}

We can now treat every surface whose anticanonical system moves: once the case \(P_B^2>0\) has been treated, \(|-K_S|\) must be a basepoint-free pencil.

\begin{corollary}\label{cor:moving-anti}
In the cyclic-cover setup, \(L\) satisfies \(\NP_p\) if \(\dim|-K_S|\geq1\).
\end{corollary}

\begin{proof}
By \cref{lem:regularity,prop:regular-case}, we may assume that \(P_B^2=0\); then \(P_B\cdot K_S=0\) by \cref{lem:isotropic-fibration}. Write \(|-K_S|=Z+|M|\), where \(Z\) is the fixed part. Since \(P_B\) is nef, \(P_B\cdot M=P_B\cdot Z=0\). The moving part \(M\) is nonzero and nef, so the Hodge index theorem gives \(M^2=0\) and \(M\equiv cP_B\) for some \(c>0\); consequently \(M\cdot Z=0\). A general member of \(|M|\) has no component in common with \(Z\), and is therefore disjoint from \(Z\). If \(Z\neq0\), their union would be a disconnected anticanonical divisor, contrary to \cref{lem:anti-subcurves}. Thus \(Z=0\). Two general members of \(|-K_S|\) have no common component and intersection number zero, so \(|-K_S|\) is basepoint-free, and it defines a morphism onto a curve. Since its members are connected, by \cref{lem:anti-subcurves}, the image is a line, so \(|-K_S|\) is a pencil. Therefore \cref{thm:elliptic-complete} applies.
\end{proof}

\section{Halphen surfaces of arbitrary index}\label{sec:halphen}

In this section \(S\) is a Halphen surface of index \(\iota\geq1\) (see \cref{sec:elliptic}), and \(F=-K_S\). The case \(\iota=1\) was treated in \cref{sec:elliptic}. When \(\iota\geq2\), the anticanonical system consists of a single curve \(F\), and the fibers of the genus-one fibration \(f\colon S\to\prj{1}\) are the members of \(|\iota F|\); thus \(\iota F\) is a multiple fiber of \(f\) with support \(F\). By Kodaira's classification of multiple fibers, \(F\) is a curve of type \(\mathrm I_N\) for some \(N\geq0\): a smooth elliptic curve (\(N=0\)), an irreducible rational curve with one node (\(N=1\)), or a cycle of \(N\geq2\) smooth rational curves; see \cite[Proposition~2.2 and Section~2.2]{CantatDolgachev}.

Since \(F\) does not move, the argument of \cref{sec:add-fibers} is replaced by adding copies of the fixed curve \(F\) one at a time (\cref{sec:add-F}). As \(F\) may be singular, the semistability argument on a smooth fiber is replaced by the elementary vanishing lemmas of \cref{sec:genus-one-curves}. Our aim is the following two theorems.

\begin{theorem}\label{thm:halphen}
In the cyclic-cover setup, if \(S\) is a Halphen surface, then \(L\) satisfies \(\NP_p\).
\end{theorem}

\begin{theorem}\label{thm:nonnegative}
In the cyclic-cover setup, if \(K_S^2\geq0\), then \(L\) satisfies \(\NP_p\).
\end{theorem}

\subsection{Two vanishing lemmas on curves of arithmetic genus one}\label{sec:genus-one-curves}

For a globally generated line bundle \(J\) on a projective curve \(C\), we write
\[
 0\longrightarrow M_J\longrightarrow H^0(C,J)\otimes\sshf{C}\longrightarrow J\longrightarrow0
\]
for its evaluation sequence. In this subsection, tensor products of line bundles are written multiplicatively. The two lemmas below treat curves of type \(\mathrm I_N\) with \(N\leq1\) and with \(N\geq2\), respectively. On such curves we do not use semistability; instead, we filter the kernel bundles by line bundles.

\begin{lemma}\label{lem:halphen-integral}
Let \(C\) be an integral projective Gorenstein curve with \(\pa(C)=1\), so that \(\omega_C\simeq\sshf{C}\). Let \(J_1,J_2\) be line bundles of the same degree \(e\geq3\), and let \(T\) be a line bundle of degree at least \(e\). If \(a,b\geq0\) and \(a+b\leq e-2\), then
\[
 H^1\bigl(C,\textstyle\bigwedge^aM_{J_1}\otimes\bigwedge^bM_{J_2}\otimes T\bigr)=0.
\]
\end{lemma}

\begin{proof}
By Serre duality, \(H^1(C,P)\simeq H^0(C,P^{-1})^\vee\) for every line bundle \(P\) on \(C\). As \(C\) is integral, \(H^1(C,P)=0\) if \(\deg P>0\), or if \(P\) is nontrivial of degree zero. Every line bundle of degree at least two is globally generated; for instance, apply \cite[Theorem 1.1]{CFHR} to clusters of length one.

For \(i=1,2\), choose \(e-2\) distinct smooth points \(x_{i,1},\ldots,x_{i,e-2}\) of \(C\), put \(Z_i=\sum_kx_{i,k}\), and let \(P_i=J_i(-Z_i)\). This is a globally generated line bundle of degree two with \(h^0(C,P_i)=2\), so \(M_{P_i}\simeq P_i^{-1}\). Since \(H^1(C,P_i)=0\), the restriction map \(H^0(C,J_i)\to H^0(J_i|_{Z_i})\) is surjective, with kernel \(H^0(C,P_i)\). Choose a basis of \(H^0(C,J_i)/H^0(C,P_i)\) that maps to a basis of \(H^0(J_i|_{Z_i})=\bigoplus_kJ_i\otimes k(x_{i,k})\) compatible with this decomposition. Comparing evaluation maps, the snake lemma gives an exact sequence
\begin{equation}\label{eq:halphen-point-filtration}
 0\longrightarrow P_i^{-1}\longrightarrow M_{J_i}
 \longrightarrow\bigoplus_{k=1}^{e-2}\sshf{C}(-x_{i,k})\longrightarrow0.
\end{equation}
Consequently, \(\bigwedge^aM_{J_i}\) has a filtration whose graded pieces are the line bundles
\[
 \sshf{C}\Bigl(-\sum_{k\in I}x_{i,k}\Bigr)\quad(|I|=a),\qquad
 P_i^{-1}\Bigl(-\sum_{k\in I}x_{i,k}\Bigr)\quad(|I|=a-1),
\]
of degrees \(-a\) and \(-a-1\), respectively; for \(a=0\), only \(\sshf{C}\) occurs. Hence every graded piece of the induced filtration of \(\bigwedge^aM_{J_1}\otimes\bigwedge^bM_{J_2}\otimes T\) is a line bundle of degree at least \(\deg T-a-b-2\geq0\), and it has no \(H^1\) unless it is trivial.

A graded piece of degree zero occurs only if \(\deg T=e\), \(a+b=e-2\), \(a,b\geq1\), and both factors come from the subbundles \(P_1^{-1}\) and \(P_2^{-1}\) in \eqref{eq:halphen-point-filtration}. Such a piece is
\[
 T\,J_1^{-1}J_2^{-1}\Bigl(Z_1+Z_2-\sum_{k\in I}x_{1,k}-\sum_{l\in I'}x_{2,l}\Bigr),
 \qquad |I|=a-1,\quad |I'|=b-1.
\]
We show that the points can be chosen so that none of these finitely many line bundles is trivial. Since \(|I|=a-1<e-2\), each of them has the form \(\Lambda(x)\), where \(x=x_{1,k}\) for some \(k\notin I\), and \(\Lambda\) does not depend on \(x\). The map \(x\mapsto\sshf{C}(x)\) from the smooth locus of \(C\) to \(\Pic^1(C)\) is injective: if \(\sshf{C}(x)\simeq\sshf{C}(y)\) with \(x\neq y\), this line bundle of degree one would have two independent sections, whereas Riemann--Roch gives \(h^0=1\). Thus, once the other points are fixed, \(\Lambda(x)\) is trivial for at most one choice of \(x\), and the configurations for which some piece is trivial form a finite union of proper closed subsets. Choosing the points outside this union, all graded pieces have vanishing \(H^1\), and the long exact sequences of the filtrations prove the lemma.
\end{proof}

\begin{lemma}\label{lem:halphen-cycle}
Let \(C=C_1+\cdots+C_N\), \(N\geq2\), be a cycle of smooth rational curves, that is, a curve of type \(\mathrm I_N\). Let \(J_1,J_2,T\) be line bundles on \(C\) such that \(J_1\) and \(J_2\) have positive degree on every component. If \(a,b\geq0\) and \(\deg(T|_{C_i})\geq a+b+1\) for every \(i\), then
\[
 H^1\bigl(C,\textstyle\bigwedge^aM_{J_1}\otimes\bigwedge^bM_{J_2}\otimes T\bigr)=0.
\]
\end{lemma}

\begin{proof}
Let \(J\) be a line bundle of positive degree \(d_i\) on each \(C_i\). We first show that the restriction maps \(H^0(C,J)\to H^0(C_i,J|_{C_i})\) are surjective; in particular, \(J\) is globally generated. Let \(C_i'\) be the union of the other components, a chain of smooth rational curves, and let \(y,y'\in C_i'\) be the points lying over the two nodes of \(C\) on \(C_i\). The kernel of the restriction map is \(H^0(C_i',J|_{C_i'}(-y-y'))\). The line bundle \(J|_{C_i'}(-y-y')\) has nonnegative degree on every component of the chain \(C_i'\) if \(N\geq3\), and degree at least \(-1\) on \(C_i'\simeq\prj{1}\) if \(N=2\). In both cases it has no \(H^1\), which proves surjectivity.

As in \cref{lem:fiber-restriction}\textup{(i)}, we obtain an exact sequence \(0\to H^0(C,J\otimes\mathscr I_{C_i})\otimes\sshf{C_i}\to M_J|_{C_i}\to M_{J|_{C_i}}\to0\), where \(\mathscr I_{C_i}\) is the ideal sheaf of \(C_i\) in \(C\) and \(M_{J|_{C_i}}\simeq\sshf{\prj{1}}(-1)^{\oplus d_i}\). It splits, as \(H^1(\prj{1},\sshf{\prj{1}}(1))=0\), so
\[
 M_J|_{C_i}\simeq\sshf{\prj{1}}^{\oplus(h^0(C,J)-d_i-1)}\oplus\sshf{\prj{1}}(-1)^{\oplus d_i}.
\]
Therefore the restriction of \(\mathcal W=\bigwedge^aM_{J_1}\otimes\bigwedge^bM_{J_2}\otimes T\) to each \(C_i\) is a direct sum of line bundles of degree at least \(\deg(T|_{C_i})-a-b\geq1\). Let \(\nu\colon\widetilde C=\bigsqcup_iC_i\to C\) be the normalization. In the cohomology sequence of
\[
 0\longrightarrow\mathcal W\longrightarrow\nu_*\nu^*\mathcal W\longrightarrow\bigoplus_{z}\mathcal W(z)\longrightarrow0,
\]
where \(z\) runs over the nodes of \(C\), we have \(H^1(\widetilde C,\nu^*\mathcal W)=0\), and the map \(H^0(\widetilde C,\nu^*\mathcal W)\to\bigoplus_z\mathcal W(z)\) is surjective: each \(C_i\) contains exactly two points lying over nodes, and a direct sum of line bundles of degree at least one on \(\prj{1}\) has sections with arbitrary values at two given points. Hence \(H^1(C,\mathcal W)=0\).
\end{proof}

\subsection{Adding copies of the anticanonical curve}\label{sec:add-F}

In this subsection, \(S\) is a smooth rational surface, and \(F\in|-K_S|\) is a nef divisor with \(F^2=0\) that is a curve of type \(\mathrm I_N\) for some \(N\geq0\). By \cref{lem:anti-subcurves}, \(F\) is Gorenstein with \(\omega_F\simeq\sshf{F}\) and \(\pa(F)=1\); and \(F\cdot C=0\) for every component \(C\) of \(F\), since \(F\) is nef and \(F^2=0\). Let \(T_0\) be one of the following:
\begin{enumerate}
\item \(T_0=0\);
\item a nonzero Cartier divisor with
\begin{equation}\label{eq:halphen-small-hypotheses}
 T_0\sim_{\mathbb{Q}}\beta F+\tfrac12R,\qquad 0\leq\beta<1,
 \quad R=\sum E_i,
\end{equation}
where the \(E_i\) are disjoint \((-2)\)-curves with \(F\cdot E_i=0\), and
\begin{equation}\label{eq:halphen-small-cohomology}
 H^0(S,\sshf{S}(K_S+T_0))=0,
 \qquad h^1(S,\sshf{S}(T_0))\leq1.
\end{equation}
\end{enumerate}
In the second case, \(H^0(S,\sshf{S}(-T_0))=H^0(S,\sshf{S}(K_S-T_0))=0\), since \(T_0\) and \(F+T_0\) are numerically nonzero effective \(\mathbb Q\)-divisors.

The next proposition plays the role of \cref{lem:add-anticanonical} when \(F\) does not move. It is proved by adding one copy of \(F\) at a time. Since \(F\) is a fixed curve that may be singular, the cohomology on \(F\) is computed with \cref{lem:halphen-integral,lem:halphen-cycle}. The induction lowers the exterior degree and at the same time replaces \(A\) by \(A-F\); for this reason, the statement is made for all \(p\geq-1\) and all \(A\) at once.

\begin{proposition}[Adding copies of \(F\)]\label{prop:halphen-induction}
Let \(p\geq-1\), and suppose that \(A\cdot C\geq p+3\) for every irreducible curve \(C\subset S\). Put
\[
 G=K_S+A=A-F,\qquad D^{(n)}=G+T_0+nF\quad(n\geq0).
\]
Then
\begin{equation}\label{eq:halphen-inductive-vanishing}
 H^1\Bigl(S,\textstyle\bigwedge^h M_{D^{(n)}}\otimes\sshf{S}(G+qD^{(n)}+tF)\Bigr)=0
\end{equation}
for all \(n,q,t\geq0\) and \(0\leq h\leq p+1\).
\end{proposition}

\begin{proof}
Put \(Q=A-(p+3)F\) and \(e=A\cdot F\geq p+3\).

\smallskip\noindent\emph{Step 0: line bundles.} The cone-theorem argument of \cref{prop:volume} shows that \(Q\) is nef. So is \(Q+T_0\): its degree on \(E_i\) is \(A\cdot E_i-1\geq p+2\geq1\), and its degree on any other curve is at least that of \(Q\), by \eqref{eq:halphen-small-hypotheses}. Hence the divisors
\[
 G=Q+(p+2)F,\qquad D^{(n)}=Q+T_0+(p+n+2)F,\qquad G+qD^{(n)}+tF
\]
are positive on every curve and have positive square, so they are ample. Their anticanonical degree is at least \(e\geq2\), so by Harbourne's theorem they are globally generated and have vanishing \(H^1\). This is \eqref{eq:halphen-inductive-vanishing} for \(h=0\), and proves the proposition when \(p=-1\).

Assume from now on that \(p\geq0\). Then \(A'=A-F\) satisfies \(A'\cdot C\geq p+2\) for every curve \(C\): if \(F\cdot C=0\), then \(A'\cdot C=A\cdot C\), and otherwise use \(A'=Q+(p+2)F\). Applying the previous paragraph to \(A'\) and \(p-1\) shows that \(D^{(n)}-F=(K_S+A')+T_0+nF\) is globally generated with \(H^1(S,\sshf{S}(D^{(n)}-F))=0\). Hence the restriction map \(H^0(S,\sshf{S}(D^{(n)}))\to H^0(F,\sshf{F}(D^{(n)}))\) is surjective, with kernel \(H^0(S,\sshf{S}(D^{(n-1)}))\) when \(n\geq1\), and as in \cref{lem:fiber-restriction}\textup{(i)} there is an exact sequence
\begin{equation}\label{eq:halphen-restriction-kernel}
 0\longrightarrow H^0(S,\sshf{S}(D^{(n)}-F))\otimes\sshf{F}
 \longrightarrow M_{D^{(n)}}|_F\longrightarrow M_{D^{(n)}|_F}\longrightarrow0.
\end{equation}
The line bundle \(\sshf{F}(D^{(n)})\) has degree \(e\geq3\), and positive degree on every component of \(F\) because \(D^{(n)}\) is ample. Moreover, for all \(q,t\geq0\), the line bundle \(\sshf{F}(G+qD^{(n)}+tF)\) has degree at least \(e\), and degree at least \(G\cdot C=A\cdot C\geq p+3\) on every component \(C\) of \(F\).

\smallskip\noindent\emph{Step 1: the case \(n=t=0\).} We show that \(H^1\bigl(S,\bigwedge^hM_{D^{(0)}}\otimes\sshf{S}(G+qD^{(0)})\bigr)=0\) for \(1\leq h\leq p+1\) and \(q\geq0\). If \(T_0=0\), then \(D^{(0)}=G\) satisfies \(\NP_p\) by \cite[Theorem 1.3]{GP01}, since \(G\cdot(-K_S)=e\geq p+3\), and the claim is the kernel-bundle form of Property \(\NP_p\). If \(T_0\neq0\), apply \cref{lem:weight-one-regularity} with \((D,T,m)=(D^{(0)},T_0,p)\), so that \(D-T=G\). Its hypotheses hold: \(D^{(0)}\) is very ample by \cite[Theorem 1.3]{GP01}, since \(D^{(0)}\cdot F=e\geq3\); \(G=K_S+A\) is \(p\)-very ample by \cref{lem:adjoint-separation}\textup{(i)}; \(D^{(0)}\cdot T_0>0\) by \eqref{eq:halphen-small-hypotheses}; the groups \(H^0(S,\sshf{S}(-T_0))\), \(H^0(S,\sshf{S}(K_S-T_0))\) and \(H^0(S,\sshf{S}(K_S+T_0))\) vanish, by the remark above and \eqref{eq:halphen-small-cohomology}; \(H^1(S,\sshf{S}(qD^{(0)}-T_0))=H^1(S,\sshf{S}(G+(q-1)D^{(0)}))=0\) for \(q\geq1\), by Step~0; and \((-K_S)\cdot G-h^1(S,\sshf{S}(T_0))\geq e-1\geq p+2\).

\smallskip\noindent\emph{Step 2: the case \(n=0\), \(t\geq0\).} We pass from \(t-1\) to \(t\) using the exact sequence
\[
 0\to\textstyle\bigwedge^hM_{D^{(0)}}\otimes\sshf{S}(\Xi-F)\to\bigwedge^hM_{D^{(0)}}\otimes\sshf{S}(\Xi)\to\bigwedge^hM_{D^{(0)}}|_F\otimes\sshf{F}(\Xi)\to0,
 \]
where \(\Xi=G+qD^{(0)}+tF\).
By \eqref{eq:halphen-restriction-kernel}, \(\bigwedge^hM_{D^{(0)}}|_F\) is filtered by the bundles \(\bigwedge^wM_{D^{(0)}|_F}\), \(w\leq h\), tensored with trivial bundles. If \(F\) is integral, \cref{lem:halphen-integral}, with \(b=0\), shows that the right-hand term has no \(H^1\), since \(w\leq h\leq p+1\leq e-2\). If \(F\) is a cycle, \cref{lem:halphen-cycle} applies, since \(\sshf{F}(\Xi)\) has degree at least \(p+3\geq w+1\) on every component.

\smallskip\noindent\emph{Step 3: induction.} We prove \eqref{eq:halphen-inductive-vanishing} for \(h\geq1\) and \(n\geq1\) by induction on \(h\) and, for fixed \(h\), on \(n\), simultaneously for all \(p\geq h-1\) and all \(A\) as in the statement. Put \(W_n=H^0(F,\sshf{F}(D^{(n)}))\), and let \(i\colon F\hookrightarrow S\) be the inclusion. Since \(D^{(n)}=D^{(n-1)}+F\), the argument of Step~2 in the proof of \cref{lem:add-anticanonical} gives an exact sequence
\[
 0\longrightarrow M_{D^{(n-1)}}\longrightarrow M_{D^{(n)}}\longrightarrow\shf{V}_n\longrightarrow0,
 \qquad
 \shf{V}_n=\ker\bigl(W_n\otimes\sshf{S}\longrightarrow i_*\sshf{F}(D^{(n)})\bigr),
\]
where \(\shf{V}_n\) is locally free and, for \(v\geq1\),
\begin{equation}\label{eq:halphen-elementary-wedge}
 0\longrightarrow\bigwedge^vW_n\otimes\sshf{S}(-F)
 \longrightarrow\bigwedge^v\shf{V}_n
 \longrightarrow i_*\bigwedge^vM_{D^{(n)}|_F}\longrightarrow0
\end{equation}
is exact. The local computation given there uses only that \(F\) is a Cartier divisor, so it applies when \(F\) is singular.

Hence \(\bigwedge^hM_{D^{(n)}}\) has a filtration with graded pieces \(\bigwedge^aM_{D^{(n-1)}}\otimes\bigwedge^v\shf{V}_n\), \(a+v=h\). Twist by \(\Xi=G+qD^{(n)}+tF\). There are three kinds of terms.
\begin{enumerate}
\item[(a)] If \(v=0\), the piece is \(\bigwedge^hM_{D^{(n-1)}}\otimes\sshf{S}(G+qD^{(n-1)}+(t+q)F)\), since \(qD^{(n)}=qD^{(n-1)}+qF\). It has no \(H^1\), by induction on \(n\).
\item[(b)] If \(v\geq1\), the subsheaf in \eqref{eq:halphen-elementary-wedge} contributes, up to a vector-space factor, the bundle \(\bigwedge^aM_{D^{(n-1)}}\otimes\sshf{S}(\Xi-F)\). Put \(A'=A-F\) and \(G'=K_S+A'=G-F\). As shown in Step~0, \(A'\) satisfies the hypothesis of the proposition with \(p'=p-1\); the corresponding divisors \(G'+T_0+nF\) are exactly \(D^{(n-1)}\); and \(\Xi-F=G'+qD^{(n-1)}+(t+q)F\). Since \(a\leq h-1\leq p'+1\), this bundle has no \(H^1\), by induction on \(h\) (for \(a=0\), by Step~0).
\item[(c)] If \(v\geq1\), the quotient in \eqref{eq:halphen-elementary-wedge} contributes the sheaf \(i_*\bigl(\bigwedge^aM_{D^{(n-1)}}|_F\otimes\bigwedge^vM_{D^{(n)}|_F}\otimes\sshf{F}(\Xi)\bigr)\); the sequence remains exact because \(\bigwedge^aM_{D^{(n-1)}}\) is locally free. By \eqref{eq:halphen-restriction-kernel} for \(D^{(n-1)}\), this bundle on \(F\) is filtered by the bundles \(\bigwedge^wM_{D^{(n-1)}|_F}\otimes\bigwedge^vM_{D^{(n)}|_F}\otimes\sshf{F}(\Xi)\), \(0\leq w\leq a\), tensored with trivial bundles. Here \(w+v\leq h\leq e-2\); the line bundles \(\sshf{F}(D^{(n-1)})\) and \(\sshf{F}(D^{(n)})\) both have degree \(e\) and are positive on every component; and \(\sshf{F}(\Xi)\) has degree at least \(e\), and at least \(p+3\geq w+v+1\) on every component. So \cref{lem:halphen-integral} or \cref{lem:halphen-cycle} shows that these bundles have no \(H^1\).
\end{enumerate}
The long exact cohomology sequences of these filtrations give \eqref{eq:halphen-inductive-vanishing} for \((h,n)\). At each step, one of the tensor factors is locally free, so we never tensor two sheaves supported on \(F\), and all the sequences used are exact.
\end{proof}

\subsection{Non-nef covering data}

\begin{lemma}\label{lem:halphen-small-data}
Let \(S\) be a Halphen surface of index \(\iota\geq2\), and let \(F\) be its unique anticanonical curve. In the cyclic-cover setup, assume that \(P_B^2=0\) and that \(B\) is not nef. Then the conclusions of \cref{lem:small-vertical} hold. In particular, \(B_0\) satisfies \eqref{eq:halphen-small-hypotheses} and \eqref{eq:halphen-small-cohomology}.
\end{lemma}

\begin{proof}
The first three paragraphs of the proof of \cref{lem:small-vertical} use only that \(F\in|-K_S|\) is nef with \(F^2=0\), that \(S\) is a rational surface containing a \((-1)\)-curve, and that \(X\) is not ruled. They therefore apply verbatim, and give \(r=2\), the stated form of \(B\) with \(s\in\{4,8\}\) and \(2b\in\mathbb Z_{>0}\), and the vanishing of \(H^0\) of \(-B_0\), \(K_S-B_0\) and \(K_S+B_0\). In particular, \(B_0\) satisfies \eqref{eq:halphen-small-hypotheses} with \(\beta=b_0\). It remains to show that \(h^1(S,\sshf{S}(B_0))\leq1\). In \cref{lem:small-vertical}, we restricted \(\sshf{S}(B_0)\) to a smooth fiber; now we must restrict to the fixed curve \(F\), which may be singular.

Suppose that \(s=4\). As in the proof of \cref{lem:small-vertical}, \(H^1(S,\sshf{S}(-B_0))=0\), and restriction to \(F\) gives
\[
 h^1(S,\sshf{S}(B_0))=h^1(F,\sshf{F}(B_0))=h^0(F,\sshf{F}(-B_0)),
\]
the last equality by Serre duality on \(F\). If \(F\) is integral, this is at most one, because \(\sshf{F}(-B_0)\) has degree zero. Suppose that \(F\) is a cycle. If no \(E_i\) is a component of \(F\), then \(R\cap F=\emptyset\), since \(E_i\cdot F=0\); so \(\sshf{F}(-B_0)\) has degree zero on every component, and again has at most one independent section. Otherwise, \(\sshf{F}(-B_0)\) has degree one on each component of \(F\) that is an \(E_i\), and degree \(-\tfrac12R\cdot C\in\{0,-1\}\) on every other component \(C\), because this degree is an integer and \(0\leq R\cdot C\leq2\). A component \(C\) adjacent to some \(E_i\subset F\) has \(R\cdot C>0\), hence degree \(-1\). A section of \(\sshf{F}(-B_0)\) therefore vanishes on every such \(C\). On each \(E_i\subset F\), it is then a section of \(\sshf{\prj{1}}(1)\) vanishing at the two points where \(E_i\) meets its neighbors, so it vanishes as well; and vanishing propagates across the remaining components, of degree zero. Hence \(h^0(F,\sshf{F}(-B_0))=0\) in this case.

Suppose that \(s=8\). As in the proof of \cref{lem:small-vertical}, an effective divisor \(V\sim B_0\) would contain \(R\), which is impossible if \(b_0=0\). If \(b_0=\tfrac12\), then \(2(V-R)+R\) would be an effective divisor in \(|F|=\{F\}\), so \(F\) would contain eight pairwise disjoint components. This is impossible. An integral \(F\) contains none of the \(E_i\). If \(F\) is a cycle of \(N\) curves, then \(N\leq9\): by Zariski's lemma, the classes of the components of the fiber \(\iota F\) are linearly independent in \(\operatorname{NS}(S)\), and they lie in \(F^\perp\), which has rank \(\rho(S)-1=9\). A cycle of \(N\leq9\) curves has at most four pairwise disjoint components. Thus \(B_0\) is not effective, and \(h^1(S,\sshf{S}(B_0))=1\) follows from Riemann--Roch as in \cref{lem:small-vertical}.
\end{proof}

\subsection{The case \texorpdfstring{\(K_S^2\geq0\)}{K\_S\^2 >= 0}}

\begin{proof}[Proof of \cref{thm:halphen}]
If \(\iota=1\), this is \cref{thm:elliptic-complete}, and if \(P_B^2>0\), it follows from \cref{lem:regularity,prop:regular-case}. Assume therefore that \(\iota\geq2\) and \(P_B^2=0\). Since \(\rho(S)=10\), \cref{lem:invariant-nonnef} shows that \(D\) satisfies \(\NP_p\), and \cref{prop:high-weights} gives condition \textup{(iii)} of \cref{thm:transfer}. By \cref{lem:kernel-criterion}, it remains to show that, for every \(1\leq j<r\),
\begin{equation}\label{eq:halphen-target}
 H^1\Bigl(S,\textstyle\bigwedge^hM_D\otimes\sshf{S}(D_j+qD)\Bigr)=0
 \qquad(1\leq h\leq p+1,\ q\geq0).
\end{equation}
In both cases below, this is \eqref{eq:halphen-inductive-vanishing} with \(t=0\), for suitable \(A\), \(T_0\) and \(n\) with \(D^{(n)}=D\).

Suppose first that \(B\) is nef. Then \(P_B=B\), and as in the proof of \cref{prop:elliptic-nef}, \cref{lem:isotropic-fibration}, the Hodge index theorem, and intersection with a \((-1)\)-curve give \(B\sim bF\) for a positive integer \(b\). For \(1\leq j<r\), put \(A_j=A+(r-1-j)bF\). Since \(F\) is nef, \(A_j\) satisfies the curve bound \(p+3\), and \(K_S+A_j=D_j\). Apply \cref{prop:halphen-induction} to \(A_j\), with \(T_0=0\) and \(n=jb\). Then \(D^{(n)}=D_j+jbF=D\), and \eqref{eq:halphen-inductive-vanishing} with \(t=0\) is \eqref{eq:halphen-target}.

If \(B\) is not nef, then \cref{lem:halphen-small-data} gives \(r=2\) and \(B=B_0+nF\), and \(D_1=K_S+A\). Apply \cref{prop:halphen-induction} to \(A\), with \(T_0=B_0\). Then \(D^{(n)}=K_S+A+B_0+nF=D\), and \eqref{eq:halphen-inductive-vanishing} with \(t=0\) is \eqref{eq:halphen-target} for \(j=1\).
\end{proof}

\begin{proof}[Proof of \cref{thm:nonnegative}]
By \cref{lem:regularity,prop:regular-case}, we may assume that \(P_B^2=0\). By \cref{lem:isotropic-fibration}, \(P_B\) defines a genus-one fibration \(f\colon S\to\prj{1}\), and \(K_S\cdot P_B=0\). By the Hodge index theorem, the intersection form is negative semidefinite on \(P_B^\perp\), with kernel spanned by \(P_B\). As \(K_S\in P_B^\perp\), it follows that \(K_S^2\leq0\), with equality only if \(K_S\) is numerically a multiple of \(P_B\). Under the hypothesis \(K_S^2\geq0\), therefore, \(K_S^2=0\) and \(-K_S\equiv\lambda P_B\), with \(\lambda>0\) because \(-K_S\) is effective and nonzero. Every \((-1)\)-curve \(E\) satisfies \(P_B\cdot E=\lambda^{-1}(-K_S\cdot E)=\lambda^{-1}>0\), so no \((-1)\)-curve is contained in a fiber of \(f\). Thus \(f\) is relatively minimal, \(S\) is a Halphen surface, and \cref{thm:halphen} applies.
\end{proof}

\section{Blow-ups of Halphen surfaces}\label{sec:fixed}

In this section we treat the remaining case
\begin{equation}\label{eq:fixed-boundary}
 K_S^2<0,\qquad h^0(S,\sshf{S}(-K_S))=1,\qquad P_B^2=0.
\end{equation}
Throughout, we work in the cyclic-cover setup and assume \eqref{eq:fixed-boundary}. We write \(F\) for the unique anticanonical divisor, which may be reducible or nonreduced, and put \(n=-K_S^2\geq1\). By \cref{lem:isotropic-fibration}, \(P_B\) defines a genus-one fibration \(f\colon S\to\prj{1}\), and every component of \(\Gamma\) is vertical. The divisor \(B\) need not be nef. Our aim is the following theorem.

\begin{theorem}
\label{thm:fixed-uniform}
In the cyclic-cover setup, assume \eqref{eq:fixed-boundary}.
Then \(L\) satisfies \(\NP_p\).
\end{theorem}

The idea is to subtract a horizontal smooth rational curve \(E\) from the character adjoints \(D_j\), and to add it back at the end. In \cref{sec:fixed-model} we describe \(S\) as a blow-up of a Halphen surface and construct \(E\). In \cref{sec:fixed-separation} we show that the divisors \(D_j-E\), \(1\leq j<r\), are nef, big and \(p\)-very ample. The curve \(E\) also supplies the bigness that \(P_B\) lacks: Kawamata--Viehweg vanishing applies to \(K_S+jB+E\), and restriction to \(F\) then gives the specialty bound required by the weight-one criterion (\cref{sec:fixed-specialty}). Finally, in \cref{sec:fixed-proof}, we add \(E\) back by restricting to \(E\simeq\prj{1}\).

\subsection{The relatively minimal model and a horizontal rational curve}\label{sec:fixed-model}

We first describe the relatively minimal model of \(f\) and the branch divisor. In this section the curves on which \(B\) is negative are denoted \(C_\alpha\) rather than \(E_\alpha\) as in \cref{prop:fractional-negative}, since the letter \(E\) is reserved for the curve of \cref{lem:fixed-threshold}.

\begin{lemma}
\label{lem:fixed-geometry}
Let \(\beta\colon S\to Y\) be the relatively minimal model of \(f\) over \(\prj{1}\), obtained by successively contracting \((-1)\)-curves contained in fibers. Then \(Y\) is a Halphen surface, say of index \(\iota\), and \(\beta\) is a composition of \(n\) blow-ups of points, each lying on the image of \(F\). Put
\[
 H=\beta^*(-K_Y),\qquad \Phi=\iota H,
\]
so that \(\Phi\) is the class of a fiber of \(f\).
\begin{enumerate}
\item The divisor \(H\) is integral and nef, \(H^2=F\cdot H=0\), and an irreducible curve \(C\) is horizontal if and only if \(H\cdot C>0\).
\item One has
\begin{equation}\label{eq:fixed-branch-form}
 \Gamma=\Phi_1+\cdots+\Phi_d+R,\qquad
 P_B\sim_{\mathbb Q}\frac d r\Phi,\qquad d\geq1,
\end{equation}
where \(\Phi_1,\dots,\Phi_d\) are distinct smooth fibers of \(f\), and \(R=\sum_\alpha C_\alpha\) is the sum of the components of \(\Gamma\) with negative self-intersection. These are the curves \(C\) with \(B\cdot C<0\); they are pairwise disjoint, and \(C_\alpha^2=-rn_\alpha\), where \(n_\alpha=-B\cdot C_\alpha\geq1\). If \(r\geq3\), then \(R\leq F\).
\item Every component of \(F\) has negative self-intersection, and
\begin{equation}\label{eq:fixed-extra-degree}
 F\cdot D_j\geq p+4\qquad(0\leq j<r).
\end{equation}
\end{enumerate}
\end{lemma}

\begin{proof}
Contracting \((-1)\)-curves contained in fibers, one at a time, gives a relatively minimal genus-one fibration on a smooth rational surface \(Y\); thus \(Y\) is a Halphen surface; see \cite[Section~2.2]{CantatDolgachev}. Write \(\beta=\sigma_1\circ\cdots\circ\sigma_k\), where \(\sigma_i\colon S_i\to S_{i-1}\) is the blow-up of a point \(x_i\), with exceptional curve \(G_i\), and \(S_0=Y\), \(S_k=S\). The image \(F_i\) of \(F\) in \(S_i\) is an effective anticanonical divisor, and
\[
 F_i=\sigma_i^*F_{i-1}-G_i.
\]
Since \(F_i\) is effective, \(x_i\in F_{i-1}\). As \(K_Y^2=0\) and \(\rho(Y)=10\), the blow-up formulas give \(k=-K_S^2=n\) and \(\rho(S)=10+n\).

(i) The class \(-K_Y\) is nef with \((-K_Y)^2=0\), and \(\iota(-K_Y)\) is the class of a fiber of \(Y\to\prj{1}\). Hence \(H\) is integral and nef, \(H^2=0\), \(F\cdot H=\beta_*F\cdot(-K_Y)=K_Y^2=0\), and \(\Phi=\iota H\) is the class of a fiber of \(f\). In particular, \(C\) is horizontal if and only if \(H\cdot C>0\).

(ii) All the points \(x_i\) lie over the fiber of \(Y\) containing \(F_0\). So the fiber \(\Phi_0\) of \(f\) containing \(F\) is reducible, and \(\beta\) is an isomorphism over the complement of the point below \(\Phi_0\). Let \(C\) be a component of \(\Gamma\). It is vertical, so \(C^2\leq0\), with equality if and only if \(C\) is the support of a fiber, by Zariski's lemma. Suppose that \(C^2=0\). Then \(C\not\subset\Phi_0\), because the components of a reducible fiber have negative square. Every multiple fiber of \(Y\to\prj{1}\) is supported on \(F_0\). Hence \(C\) is a fiber of \(f\), and it is smooth, since it is a component of the smooth divisor \(\Gamma\). Suppose that \(C^2<0\). As \(\Gamma\) is smooth, \(C^2=\Gamma\cdot C=rB\cdot C\), so \(B\cdot C<0\). Conversely, every curve on which \(B\) is negative is such a component, by \cref{prop:fractional-negative}. This gives the decomposition of \(\Gamma\), and the assertions about the \(C_\alpha\) are \cref{prop:fractional-negative}. By \eqref{eq:zariski-B}, \(P_B=B-\frac1rR\sim_{\mathbb Q}\frac1r(\Phi_1+\cdots+\Phi_d)\), and \(d\geq1\) because \(P_B\not\equiv0\).

If some \(C_\alpha\) is not contained in \(F\), then \(F\cdot C_\alpha\geq0\), and adjunction gives
\[
 0\leq F\cdot C_\alpha=2-2g(C_\alpha)-rn_\alpha.
\]
Hence \(r=2\), and \(C_\alpha\) is a \((-2)\)-curve. Thus \(R\leq F\) if \(r\geq3\).

(iii) If \(F\) is integral, then \(F^2=-n<0\). Otherwise, every component \(C\) of \(F\) is a smooth rational curve, by \cref{lem:anti-subcurves}\textup{(ii)}. If \(C^2\geq0\), then \(h^2(S,\sshf{S}(C))=h^0(S,\sshf{S}(K_S-C))=0\), so \(h^0(S,\sshf{S}(C))\geq\chi(\sshf{S}(C))=C^2+2\geq2\). Then \(F=C+(F-C)\) would move, contrary to \(h^0(S,\sshf{S}(-K_S))=1\).

Since \(\rho(S)\geq3\), \cref{lem:invariant-nonnef} applies. If \(F\) is not integral, it gives \(F\cdot D_j\geq2(p+2)\geq p+4\). Suppose that \(F\) is integral, and put \(t=r-1-j\). If \(B\cdot F\geq0\), then
\(F\cdot D_j=-F^2+A\cdot F+tB\cdot F\geq n+p+3\geq p+4\).
If \(B\cdot F=-a<0\), then \(F\) is a component of \(\Gamma\) with \(F^2=-ra\), and
\(F\cdot D_j=A\cdot F+(j+1)a\geq p+4\).
\end{proof}

We next lower the polarization \(A\) in the fiber direction until the bound \(A\cdot C\geq p+3\) is attained on a horizontal curve. This curve will be subtracted from the character adjoints.

\begin{lemma}
\label{lem:fixed-threshold}
There are a rational number \(\tau\geq0\) and a smooth rational curve \(E\), horizontal for \(f\), such that the class \(J=A-\tau H\) satisfies
\begin{equation*}
 J\cdot C\geq p+3\ \text{for every curve }C,\qquad
 J\cdot E=p+3,\qquad E^2\in\{-1,-2\}.
\end{equation*}
Moreover, \(J\) is ample, \(J-(p+3)F\) is nef, and
\begin{equation*}
 J^2\geq(n+2)(p+3)^2\geq3(p+3)^2.
\end{equation*}
\end{lemma}

\begin{proof}
\emph{Step 1: real classes satisfying the curve bound.} Let \(J\) be an \(\mathbb R\)-divisor class with \(J\cdot C\geq p+3\) for every irreducible curve \(C\). Then \(J\) is nef. The cone-theorem argument in the proof of \cref{prop:volume} does not use integrality, and \(\rho(S)\geq3\); it shows that \(J-(p+3)F\) is nef. Since \(F^2=-n\) and \(J\cdot F\geq p+3\),
\begin{align*}
 J^2&=\bigl(J-(p+3)F\bigr)^2+2(p+3)J\cdot F+(p+3)^2n\\
 &\geq(n+2)(p+3)^2.
\end{align*}
As \(J\) is positive on every curve and \(J^2>0\), the Nakai--Moishezon criterion for \(\mathbb R\)-divisors \cite[Theorem~1.4]{FujinoMiyamoto} shows that \(J\) is ample.

\smallskip\noindent\emph{Step 2: the threshold.} Define
\begin{equation}\label{eq:fixed-threshold}
 \tau=\inf_C\frac{A\cdot C-p-3}{H\cdot C},\qquad J=A-\tau H,
\end{equation}
where \(C\) runs over the horizontal irreducible curves. Since \(A\cdot C\geq p+3\) for every curve, and horizontal curves exist, \(0\leq\tau<\infty\). By definition, \(J\cdot C\geq p+3\) for every horizontal irreducible curve \(C\); if \(C\) is vertical, then \(J\cdot C=A\cdot C\geq p+3\). Hence Step~1 applies to \(J\).

\smallskip\noindent\emph{Step 3: the infimum is attained.} We first show that only finitely many numerical classes of irreducible curves \(C\) satisfy \(J\cdot C<2(p+3)\). Fix an integral ample divisor \(T\), and choose \(c,a>0\) such that \(cJ-T\) and \(aT\pm K_S\) are ample. For such a curve \(C\), the degree \(T\cdot C\leq cJ\cdot C\) is bounded, \(|K_S\cdot C|\leq aT\cdot C\), and adjunction gives \(C^2\geq-K_S\cdot C-2\geq-aT\cdot C-2\). Write
\[
 C\equiv\frac{T\cdot C}{T^2}\,T+C_\perp,\qquad C_\perp\in T^\perp.
\]
Then \(C_\perp^2=C^2-(T\cdot C)^2/T^2\) is bounded below. The intersection form is negative definite on \(T^\perp\), so these classes lie in a bounded subset of the N\'eron--Severi lattice, and hence form a finite set.

Suppose now that \(J\cdot C>p+3\) for every horizontal irreducible curve \(C\). By finiteness, the ratios 
\[\frac{J\cdot C-p-3}{H\cdot C},\] as \(C\) ranges over horizontal irreducible curves satisfying \(J\cdot C<2(p+3)\), have a positive minimum \(\mu\). If there is no such curve, put \(\mu=1\). Choose \(\eta>0\) such that \(J-\eta H\) is ample, and \(0<\varepsilon\leq\min\{\mu,\eta/2\}\). If \(J\cdot C\geq2(p+3)\), then, writing \(J-\varepsilon H=(1-\varepsilon/\eta)J+(\varepsilon/\eta)(J-\eta H)\), we get
\[
 (J-\varepsilon H)\cdot C\geq(1-\varepsilon/\eta)\,J\cdot C\geq J\cdot C/2\geq p+3.
\]
If \(C\) is horizontal with \(J\cdot C<2(p+3)\), then \((J-\varepsilon H)\cdot C\geq p+3\) because \(\varepsilon\leq\mu\); and degrees on vertical curves do not change. Thus \(\tau+\varepsilon\) also satisfies the bound, contrary to \eqref{eq:fixed-threshold}. Hence there is a horizontal irreducible curve \(E\) with \(J\cdot E=p+3\). Then \(\tau=(A\cdot E-p-3)/(H\cdot E)\) is rational.

\smallskip\noindent\emph{Step 4: the curve \(E\).} If \(E^2>0\), the Hodge index theorem gives \(J\cdot E\geq\sqrt{J^2E^2}\geq\sqrt3\,(p+3)>p+3\), which is impossible.

We show that \(J\cdot C>p+3\) for every irreducible curve \(C\) with \(C^2=0\); this excludes \(E^2=0\). Such a curve is not contained in \(F\), by \cref{lem:fixed-geometry}\textup{(iii)}, so \(F\cdot C\geq0\), and adjunction gives \(F\cdot C=2-2\pa(C)\). If \(\pa(C)=0\), then \(J\cdot C\geq(p+3)F\cdot C=2(p+3)\), as \(J-(p+3)F\) is nef. If \(\pa(C)=1\), then \(F\cdot C=0\), so \(C\) is disjoint from \(F\). The exact sequence
\[
 0\longrightarrow\sshf{S}(K_S+C)\longrightarrow\sshf{S}(C)\longrightarrow\sshf{F}\longrightarrow0,
\]
together with \(H^2(S,\sshf{S}(K_S+C))\simeq H^0(S,\sshf{S}(-C))^\vee=0\), gives \(h^1(S,\sshf{S}(C))\geq h^1(F,\sshf{F})=1\). Since \(\chi(\sshf{S}(C))=1\) and \(h^2(S,\sshf{S}(C))=0\), we get \(h^0(S,\sshf{S}(C))\geq2\). A second member of \(|C|\) has no component in common with \(C\), and is therefore disjoint from \(C\), because \(C^2=0\). These two members span a basepoint-free pencil, whose fibers are connected, since one of them is the integral curve \(C\). So the pencil defines a genus-one fibration \(f'\colon S\to\prj{1}\) with \(C\) as a fiber. Let \(\beta'\colon S\to Y'\) be its relatively minimal model, a Halphen surface of some index \(\iota'\). As in the proof of \cref{lem:fixed-geometry}, \(\beta'\) is a composition of \(n\) point blow-ups, and \(K_{S/Y'}\) is the sum of the total transforms of their exceptional curves. These are \(n\) nonzero effective divisors, so \(J\cdot K_{S/Y'}\geq n(p+3)\). Since \(C\equiv\iota'\beta'^*(-K_{Y'})=\iota'(F+K_{S/Y'})\),
\[
 J\cdot C\geq\iota'\bigl(J\cdot F+n(p+3)\bigr)\geq(n+1)(p+3)>p+3.
\]

Thus \(E^2<0\). As \(E\) is horizontal, it is not contained in \(F\), so \(F\cdot E\geq0\), and adjunction gives \(E^2=2\pa(E)-2+F\cdot E\). Hence \(\pa(E)=0\) and \(E^2\in\{-1,-2\}\); that is, \(E\) is a smooth rational curve.
\end{proof}

From now on, we fix \(\tau\), \(J\) and \(E\) as in \cref{lem:fixed-threshold}, and write
\[
 E^2=-m,\qquad m\in\{1,2\},\qquad\text{so that}\quad F\cdot E=2-m\leq1.
\]
Since \(E\) is horizontal and every component of \(F\) and of \(\Gamma\) is vertical, \(E\) is a component of neither \(F\) nor \(\Gamma\). In particular, \(E\cdot C\leq F\cdot E\leq1\) for every component \(C\) of \(F\), and \(B\cdot E\geq0\).

\subsection{Separation after subtracting the rational curve}\label{sec:fixed-separation}

We now show that subtracting \(E\) leaves enough positivity for the weight-one criterion. The class \(J\) enters only through numerical estimates; the separation and vanishing theorems are applied only to integral divisors.

\begin{proposition}
\label{prop:fixed-predecessors}
Put \(U=A-E\) and \(P=K_S+U=D_{r-1}-E\). Then \(U\) is ample, \(U^2>(p+2)^2\), and
\begin{equation}\label{eq:fixed-U-bounds}
 \begin{gathered}
 U\cdot C\geq p+1\quad\text{for every curve }C,\\
 U\cdot C\geq p+2\quad\text{for every irreducible curve \(C\) with \(C^2\geq0\) or \(C\subset F\).}
 \end{gathered}
\end{equation}
The divisors \(P_j=D_j-E\), \(1\leq j<r\), are nef, big, and \(p\)-very ample.
\end{proposition}

\begin{proof}
Put \(s=J^2\geq3(p+3)^2\) and \(U_0=J-E\), so that \(U=U_0+\tau H\).

\smallskip\noindent\emph{Step 1: the bounds \eqref{eq:fixed-U-bounds}.} Let \(C\neq E\) be an irreducible curve not contained in \(F\), and put \(x=J\cdot C\geq p+3\) and \(k=E\cdot C\geq0\). Since \(F\cdot C\geq0\), adjunction gives \(C^2\geq-2\). Write
\[
 E\equiv\frac{p+3}{s}J+E',\qquad C\equiv\frac{x}{s}J+C',\qquad E',C'\in J^\perp.
\]
By the Hodge index theorem, the intersection form is negative definite on \(J^\perp\), so
\[
 k=\frac{(p+3)x}{s}+E'\cdot C'\leq\frac{(p+3)x}{s}+\sqrt{(-E'^2)(-C'^2)},
\]
where \(-E'^2=m+(p+3)^2/s\) and \(-C'^2=-C^2+x^2/s\leq2+x^2/s\). Hence
\begin{equation*}
 k\leq\psi(x):=\frac{(p+3)x}{s}
       +\sqrt{\left(m+\frac{(p+3)^2}{s}\right)
                     \left(2+\frac{x^2}{s}\right)}.
\end{equation*}
Since \(s\geq3(p+3)^2\), \(m\leq2\), and \(\frac{d}{dx}\sqrt{2+x^2/s}\leq1/\sqrt s\), one has
\[
 \psi(p+3)\leq\frac83<3,\qquad
 \psi'(x)\leq\theta:=\frac{p+3+\sqrt{(p+3)^2+ms}}s
       \leq\frac{1+\sqrt7}{3(p+3)}<1.
\]
Therefore \(k\leq\psi(p+3)+\theta(x-p-3)<x-p\), that is, \(U_0\cdot C>p\). Since \(\tau H\) is nef and \(U\) is integral, \(U\cdot C\geq p+1\).

If moreover \(C^2\geq0\), the same computation with \(-C'^2\leq x^2/s\) gives \(k\leq\theta x\). As \(\theta(p+3)\leq(1+\sqrt7)/3<2\), we get \(x-k\geq(1-\theta)x\geq(p+3)-\theta(p+3)>p+1\), so \(U\cdot C\geq U_0\cdot C>p+1\), and \(U\cdot C\geq p+2\).

On \(E\), one has \(U_0\cdot E=p+3+m\). If \(C\) is a component of \(F\), then \(H\cdot C=0\) and \(E\cdot C\leq1\), so \(U\cdot C=U_0\cdot C=A\cdot C-E\cdot C\geq p+2\). Thus \(U_0\) is nef, \eqref{eq:fixed-U-bounds} holds, and, since \(H\) is nef with \(H^2=0\),
\[
 U^2=U_0^2+2\tau\,U_0\cdot H\geq U_0^2=s-2(p+3)-m
   \geq3(p+3)^2-2(p+3)-2>(p+2)^2.
\]
By the Nakai--Moishezon criterion, \(U\) is ample.

\smallskip\noindent\emph{Step 2: the divisor \(P\).} Put \(\ell=p+1\). We apply the Reider criterion \eqref{eq:AB-Reider} to \(U\), with \(d=\ell\). We have \(U^2>(\ell+1)^2\). Let \(C\) be a nonzero effective divisor with \(C^2\leq0\). If \(C\) is irreducible with \(C^2<0\), then \(U\cdot C\geq\ell>C^2+\ell\). If \(C\) is irreducible with \(C^2=0\), then \(U\cdot C\geq\ell+1>C^2+\ell\). Otherwise \(C\) has at least two components, counted with multiplicity, and \(U\cdot C\geq2\ell>\ell\geq C^2+\ell\). Hence \(P=K_S+U\) separates every cluster of length at most \(p+1\); that is, \(P\) is \(p\)-very ample. In particular, \(P\) is globally generated, hence nef.

To see that \(P\) is big, note that \(K_S\cdot H=-F\cdot H=0\), so \((K_S+U_0)\cdot H=U_0\cdot H\geq0\). As \(P=(K_S+U_0)+\tau H\) and \(H^2=0\), this gives \(P^2\geq(K_S+U_0)^2\). Since \(\bigl(J-(p+3)F\bigr)^2=s-2(p+3)J\cdot F-(p+3)^2n\geq0\) by \cref{lem:fixed-threshold}, and \(F\cdot E=2-m\),
\begin{align*}
 P^2\geq(K_S+U_0)^2&=s-n-m-2J\cdot F-2(p+3)+2F\cdot E\\
 &\geq\left(1-\frac1{p+3}\right)s+(p+2)n-2(p+3)+4-3m\\
 &\geq3(p+3)^2-4(p+3)-3>0.
\end{align*}

\smallskip\noindent\emph{Step 3: the divisors \(P_j\).} Write \(P_j=K_S+U+tB\), where \(t=r-1-j\), so \(0\leq t\leq r-2\). The case \(t=0\) is Step~2. If \(B\) is nef, then \(U+tB\) is ample, \((U+tB)^2\geq U^2\) and \((U+tB)\cdot C\geq U\cdot C\) for every curve \(C\), so the argument of Step~2 applies to \(U+tB\). Assume from now on that \(B\) is not nef and that \(t\geq1\). Then \(r\geq3\), so \(R\leq F\) by \cref{lem:fixed-geometry}\textup{(ii)}.

Put \(a_\alpha=U\cdot C_\alpha\) and
\[
 \eta_\alpha=\min\left\{\frac tr,\frac{a_\alpha}{rn_\alpha}\right\},
 \qquad
 P_{U+tB}=U+tP_B+\sum_\alpha\eta_\alpha C_\alpha.
\]
Since \(U+tB=U+tP_B+\frac trR\), this is the positive part of the Zariski decomposition of \(U+tB\). Indeed, \(P_{U+tB}\cdot C\geq U\cdot C\) for every irreducible curve \(C\) other than the \(C_\alpha\), and \(P_{U+tB}\cdot C_\alpha=\max\{a_\alpha-tn_\alpha,0\}\); so \(P_{U+tB}\) is nef, and the negative part \(\sum_\alpha(\frac tr-\eta_\alpha)C_\alpha\) is supported on pairwise disjoint curves on which \(P_{U+tB}\) has degree zero. Since \(P_B^2=P_B\cdot C_\alpha=0\),
\begin{equation*}
 P_{U+tB}^2=U^2+2tU\cdot P_B+
   \sum_\alpha\eta_\alpha(2a_\alpha-rn_\alpha\eta_\alpha)
   \geq U^2>(\ell+1)^2.
\end{equation*}

Suppose that \(K_S+U+tB\) fails to separate some cluster of length at most \(\ell\), and let \(Z\) be a minimal subcluster on which separation fails. By \eqref{eq:Langer-cluster-bound} and the discussion following it, \(Z\) is Gorenstein, and its Langer invariant is \(\delta=4\operatorname{length}(Z)\leq4\ell\). By a theorem of Langer \cite[Theorem~0.1]{Langer}, applied to the pseudoeffective divisor \(U+tB\), there is a nonzero effective divisor \(V\) such that
\begin{equation}\label{eq:fixed-Langer-obstruction}
 \begin{gathered}
 (U+tB-V)\cdot C\leq\delta/4\leq\ell
       \quad\text{for every nonzero subcurve }C\leq V,\\
 P_{U+tB}\cdot V\leq\frac{\delta/2}{1+\sqrt{1-\delta/P_{U+tB}^2}}<\ell+1.
 \end{gathered}
\end{equation}
For the last inequality, note that the middle term is the smaller root of \(y-y^2/P_{U+tB}^2=\delta/4\). The function \(y\mapsto y-y^2/P_{U+tB}^2\) is increasing for \(y\leq P_{U+tB}^2/2\); moreover \(\ell+1<P_{U+tB}^2/2\), and its value at \(y=\ell+1\) exceeds \(\ell\geq\delta/4\), because \((\ell+1)^2<P_{U+tB}^2\).

Write \(V=V_0+\sum_\alpha b_\alpha C_\alpha\), where no \(C_\alpha\) is a component of \(V_0\). As \(P_{U+tB}\) is nef, and \(P_{U+tB}\cdot C\geq U\cdot C\geq\ell\) for every irreducible curve \(C\) other than the \(C_\alpha\), we have \(\ell+1>P_{U+tB}\cdot V\geq U\cdot V_0\). Since \(2\ell\geq\ell+1\), it follows that \(V_0\) is either zero or an irreducible curve \(C_0\) with \(U\cdot C_0=\ell\). In the latter case, \(C_0^2<0\) and \(C_0\not\subset F\) by \eqref{eq:fixed-U-bounds}; as \(R\leq F\),
\[
 R\cdot C_0\leq F\cdot C_0=C_0^2+2-2\pa(C_0)\leq1.
\]
Suppose that \(b_\alpha\geq1\) for some \(\alpha\). Then \(C_\alpha\leq V\), and since \(t\leq r-2\), \(a_\alpha\geq\ell\) and \(V_0\cdot C_\alpha\leq1\),
\[
 (U+tB-V)\cdot C_\alpha
 =a_\alpha+(rb_\alpha-t)n_\alpha-V_0\cdot C_\alpha
 \geq\ell+2n_\alpha-1>\ell,
\]
contrary to \eqref{eq:fixed-Langer-obstruction}. Hence \(V=C_0\). But \(C_0\) is not one of the \(C_\alpha\), so \(B\cdot C_0\geq0\), and
\[
 (U+tB-C_0)\cdot C_0\geq U\cdot C_0-C_0^2>\ell,
\]
again contrary to \eqref{eq:fixed-Langer-obstruction}. Thus \(P_j\) is \(p\)-very ample, and in particular nef. Finally, \(P\cdot B\geq0\), as \(P\) is nef and \(B\sim_{\mathbb Q}\Gamma/r\) is \(\mathbb Q\)-effective. Hence \(P_j\cdot P=P^2+tP\cdot B>0\) and \(P_j^2=P_j\cdot P+tP_j\cdot B>0\), so \(P_j\) is big.
\end{proof}

\subsection{A specialty bound}\label{sec:fixed-specialty}

To apply \cref{prop:weight-one} to the divisors $P_j$, it remains to establish the specialty bound. We obtain it by restricting to the full anticanonical divisor $F$.

\begin{lemma}
\label{lem:fixed-augmented-specialty}
For \(1\leq j<r\), put \(T_j=jB+E\), so that \(P_j=D-T_j\). Then
\begin{equation*}
 h^1(S,\sshf{S}(T_j))=h^0(F,\sshf{F}(-T_j)),\qquad
 F\cdot P_j-h^1(S,\sshf{S}(T_j))\geq p+2.
\end{equation*}
\end{lemma}

\begin{proof}
\emph{Step 1: a bound for \(h^0(F,\sshf{F}(-T_j))\).} By \eqref{eq:fixed-extra-degree}, \(F\cdot P_j=F\cdot D_j-F\cdot E\geq p+3\), and by \cref{lem:invariant-nonnef}, \(P_j\cdot C=D_j\cdot C-E\cdot C\geq p+1\) for every component \(C\) of \(F\). Fix a component \(C_1\) of \(F\). Since \(D\) is very ample, a general section \(s\in H^0(S,\sshf{S}(D))\) restricts to a nonzerodivisor on \(F\), and has \(D\cdot C_1\geq p+2\) distinct zeros on \(C_1\), none of which is a singular point of the reduced curve \(F_{\mathrm{red}}\). Let \(Z\) be the reduced cluster formed by \(p+2\) of them. For every proper positive subcycle \(Y<F\), one has \(\pa(Y)\leq0\) by \cref{lem:anti-subcurves}\textup{(ii)}. If \(Y\) contains \(C_1\), then
\[
 P_j\cdot Y\geq p+1>2\pa(Y)-2+\operatorname{length}(Z\cap Y).
\]
Otherwise \(Z\cap Y=\varnothing\), and the same strict inequality holds since \(P_j\) is nef. For \(Y=F\), we have \(P_j\cdot F\geq p+3>2\pa(F)-2+\operatorname{length}(Z)\). Automatic adjunction \cite[Theorem~1.1 and Lemma~2.4]{CFHR} therefore gives
\(H^1(F,\is{Z}\otimes\sshf{F}(P_j))=0\). Hence \(H^1(F,\sshf{F}(P_j))=0\), so \(h^0(F,\sshf{F}(P_j))=F\cdot P_j\), and evaluation at \(Z\) is surjective. Multiplication by \(s\) maps \(H^0(F,\sshf{F}(-T_j))\) injectively into the kernel of this evaluation map. Therefore
\begin{equation}\label{eq:fixed-augmented-section-bound}
 h^0(F,\sshf{F}(-T_j))\leq F\cdot P_j-(p+2).
\end{equation}

\smallskip\noindent\emph{Step 2: vanishing.} We show that \(H^i(S,\sshf{S}(K_S+T_j))=0\) for \(i>0\). Granting this, the cohomology sequence of
\[
 0\longrightarrow\sshf{S}(K_S+T_j)\longrightarrow\sshf{S}(T_j)\longrightarrow\sshf{F}(T_j)\longrightarrow0
\]
and Serre duality on \(F\) give \(h^1(S,\sshf{S}(T_j))=h^1(F,\sshf{F}(T_j))=h^0(F,\sshf{F}(-T_j))\), and the lemma follows from \eqref{eq:fixed-augmented-section-bound}.

Since \(E\) is horizontal, \(P_B\cdot E>0\). Put
\[
 \lambda_j=\min\{1,jP_B\cdot E/m\}\in(0,1],\qquad
 Q_j=jP_B+\lambda_jE=jB-\frac jrR+\lambda_jE.
\]
Then \(Q_j\cdot E=jP_B\cdot E-m\lambda_j\geq0\), and \(Q_j\cdot C\geq0\) for every other irreducible curve \(C\), so \(Q_j\) is nef; and \(Q_j^2=\lambda_j(2jP_B\cdot E-m\lambda_j)>0\), so \(Q_j\) is big. Moreover,
\[
 \ceil{Q_j}=jB+E=T_j,\qquad
 \Delta_j:=T_j-Q_j=\frac jrR+(1-\lambda_j)E.
\]
Since \(K_S+T_j=K_S+\Delta_j+Q_j\), Kawamata--Viehweg vanishing for klt pairs \cite[Corollary~5.7.7]{FujinoKV} gives the claim once we know that \((S,\Delta_j)\) is klt.

The coefficients of \(\Delta_j\) are \(j/r<1\) and \(1-\lambda_j<1\). The components of \(R\) are smooth and pairwise disjoint, and \(E\) is smooth; so at most two components of \(\Delta_j\), namely \(E\) and one \(C_\alpha\), pass through any point. It therefore suffices to show that, at each point of \(E\cap R\), either the two curves meet transversally, or both coefficients are at most \(\frac12\). In the latter case the pair is klt even if the curves are tangent: on the blow-ups separating them, the exceptional curves receive coefficients at most zero.

If \(r\geq3\), then \(R\leq F\), so \(R\cdot E\leq F\cdot E\leq1\), and all intersections are transverse. Suppose that \(r=2\). Then \(j=1\), the coefficient of \(R\) is \(\frac12\), and by \eqref{eq:fixed-branch-form},
\[
 P_B\cdot E=\frac{d\iota}{2}\,H\cdot E,\qquad H\cdot E\in\mathbb Z_{>0}.
\]
If \(\lambda_1\geq\frac12\), both coefficients are at most \(\frac12\). If \(\lambda_1<\frac12\), then \(m=2\), since \(\lambda_1\geq\min\{1,\frac12\}\) when \(m=1\); hence \(d\iota H\cdot E<2\), so \(d=\iota=H\cdot E=1\) and \(E\cdot\Phi=1\). For a component \(C_\alpha\) of \(R\), let \(\Phi_y=\sum_v\mu_vC_v\) be the fiber containing it. As \(E\) is horizontal, \(E\cdot C_v\geq0\) for all \(v\), so
\[
 0\leq\mu_\alpha E\cdot C_\alpha\leq E\cdot\Phi_y=E\cdot\Phi=1.
\]
Thus \(E\cdot C_\alpha\leq1\), and \(E\) meets \(C_\alpha\) transversally.
\end{proof}

\subsection{Restriction to the rational curve}\label{sec:fixed-proof}

\begin{proof}[Proof of \cref{thm:fixed-uniform}]
Since \(\rho(S)\geq3\), \cref{lem:invariant-nonnef} shows that \(D\) satisfies \(\NP_p\); in particular, \(D\) is very ample. This is condition \textup{(i)} of \cref{thm:transfer}, and \cref{prop:high-weights} gives condition \textup{(iii)}. Put \(V=H^0(S,\sshf{S}(D))\). By the identification preceding \cref{lem:kernel-criterion}, \(K_{h-1,1}(N_{D,jB};V)\) is a subspace of \(H^1\bigl(S,\bigwedge^hM_D\otimes\sshf{S}(D_j)\bigr)\). So condition \textup{(ii)} follows once we prove that
\begin{equation}\label{eq:fixed-target}
 H^1\Bigl(S,\textstyle\bigwedge^hM_D\otimes\sshf{S}(D_j)\Bigr)=0
 \qquad(1\leq h\leq p+1,\ 1\leq j<r).
\end{equation}
Fix \(1\leq j<r\), and let \(T_j=jB+E\) and \(P_j=D_j-E=D-T_j\) be as in \cref{lem:fixed-augmented-specialty}.

\smallskip\noindent\emph{Step 1: subtracting \(E\).} We show that
\begin{equation}\label{eq:fixed-predecessor-kernel}
 H^1\Bigl(S,\textstyle\bigwedge^hM_D\otimes\sshf{S}(P_j)\Bigr)=0
 \qquad(1\leq h\leq p+1).
\end{equation}
We apply \cref{prop:weight-one} to \((D,T_j)\), so that \(G=P_j\). Its hypotheses hold. The divisor \(D\) is very ample, and \(P_j\) is nef and big, and \(p\)-very ample by \cref{prop:fixed-predecessors}. Since \(D\) is ample and \(B\sim_{\mathbb Q}\Gamma/r\), we have \(D\cdot T_j=jD\cdot B+D\cdot E>0\); and \(-T_j\) and \(K_S-T_j\) have negative degree on \(D\), so \(H^0(S,\sshf{S}(-T_j))=H^0(S,\sshf{S}(K_S-T_j))=0\). For every integer \(a\geq1\), the divisor \(aD-T_j=P_j+(a-1)D\) is nef with anticanonical degree at least \(F\cdot P_j\geq p+3\), so its \(H^1\) vanishes by Harbourne's theorem. The specialty bound is \cref{lem:fixed-augmented-specialty}. Hence \(K_{a,1}(N_{D,T_j};V)=0\) for \(0\leq a\leq p\), where \(N_{D,T_j}=\bigoplus_{q\geq0}H^0(S,\sshf{S}(P_j+qD))\). As \(H^1(S,\sshf{S}(P_j))=0\), the identification preceding \cref{lem:kernel-criterion} gives \(H^1\bigl(S,\bigwedge^hM_D\otimes\sshf{S}(P_j)\bigr)\simeq K_{h-1,1}(N_{D,T_j};V)\), which proves \eqref{eq:fixed-predecessor-kernel}.

\smallskip\noindent\emph{Step 2: the restriction of \(M_D\) to \(E\).} We shall show that \(D-E\) is nef and big, with \(H^1(S,\sshf{S}(D-E))=0\). Note that \(D-E=P+(r-1)B\). Let \(C\) be an irreducible curve. If \(B\cdot C\geq0\), then \((D-E)\cdot C\geq P\cdot C\geq0\). Otherwise \(C=C_\alpha\) for some \(\alpha\). If \(C_\alpha\not\subset F\), then \(r=2\), \(C_\alpha^2=-2\), \(F\cdot C_\alpha=0\) and \(B\cdot C_\alpha=-1\), by the proof of \cref{lem:fixed-geometry}\textup{(ii)}; hence \((D-E)\cdot C_\alpha=U\cdot C_\alpha-1\geq p\). If \(C_\alpha\subset F\), then adjunction gives \(K_S\cdot C_\alpha=2g(C_\alpha)-2+rn_\alpha\), and since \(E\cdot C_\alpha\leq1\),
\[
 (D-E)\cdot C_\alpha=A\cdot C_\alpha+2g(C_\alpha)-2+n_\alpha-E\cdot C_\alpha\geq p+n_\alpha\geq1.
\]
Thus \(D-E\) is nef. It is also big, because
\begin{align*}
 (D-E)^2&=(D-E)\cdot P+(r-1)(D-E)\cdot B\\
 &\geq (D-E)\cdot P\\
 &=P^2+(r-1)P\cdot B\\
 &\geq P^2>0.
\end{align*}
By \eqref{eq:fixed-extra-degree}, \(F\cdot(D-E)\geq p+4-1\), so \(H^1(S,\sshf{S}(D-E))=0\) by Harbourne's theorem. Consequently the restriction map \(H^0(S,\sshf{S}(D))\to H^0(E,\sshf{E}(D))\) is surjective, and, as in \cref{lem:fiber-restriction}\textup{(i)}, there is an exact sequence
\[
 0\longrightarrow H^0(S,\sshf{S}(D-E))\otimes\sshf{E}\longrightarrow M_D|_E\longrightarrow M_{D|_E}\longrightarrow0.
\]
Since \(E\simeq\prj{1}\), we have \(M_{D|_E}\simeq\sshf{\prj{1}}(-1)^{\oplus D\cdot E}\), and the sequence splits, as \(H^1(\prj{1},\sshf{\prj{1}}(1))=0\). Thus \(M_D|_E\) is a direct sum of copies of \(\sshf{\prj{1}}\) and \(\sshf{\prj{1}}(-1)\).

\smallskip\noindent\emph{Step 3: adding \(E\) back.} We have \(A\cdot E=J\cdot E+\tau H\cdot E\geq p+3\), \(K_S\cdot E=m-2\), and \(B\cdot E\geq0\). Hence
\[
 D_j\cdot E=A\cdot E+m-2+(r-1-j)B\cdot E\geq p+2.
\]
By Step~2, \(\bigl(\bigwedge^hM_D\otimes\sshf{S}(D_j)\bigr)|_E\) is a direct sum of line bundles of degree at least \(p+2-h\geq1\), so it has no \(H^1\). The cohomology sequence of
\[
 0\longrightarrow\textstyle\bigwedge^hM_D\otimes\sshf{S}(P_j)\longrightarrow\bigwedge^hM_D\otimes\sshf{S}(D_j)\longrightarrow\bigl(\bigwedge^hM_D\otimes\sshf{S}(D_j)\bigr)|_E\longrightarrow0
\]
and \eqref{eq:fixed-predecessor-kernel} prove \eqref{eq:fixed-target}.
\end{proof}

\section{Low Picard number and completion of the proof}\label{sec:completion}

It remains to treat surfaces of Picard number at most two, and to put together the preceding cases. If \(\rho(S)\leq2\), then \(S\) is \(\prj{2}\) or \(\mathbb F_e\), and on \(\mathbb F_e\) the kernel bundles can be computed explicitly.

\begin{lemma}
\label{lem:hirzebruch-kernel}
Let \(\varphi\colon S=\mathbb F_e\to\prj{1}\) be a Hirzebruch surface, with \(e\geq0\). Let \(C_0\) be a section with \(C_0^2=-e\), let \(f\) be a fiber, and put \(\xi=C_0+ef\). For integers \(u,v\geq0\), the line bundle \(\shf{G}=\sshf{S}(u\xi+vf)\) is globally generated, and its kernel bundle fits into an exact sequence
\begin{equation}\label{eq:Fe-evaluation}
 0\longrightarrow
 \sshf{S}(-f)^{\oplus N_{\shf{G}}}
 \longrightarrow M_{\shf{G}}
 \longrightarrow
 \bigoplus_{t=1}^u
 \sshf{S}\bigl(-\xi+(v+te)f\bigr)
 \longrightarrow0,
\end{equation}
where \(N_{\shf{G}}=\sum_{i=0}^u(v+ie)\). For \(u=0\), the direct sum on the right is zero.
\end{lemma}

\begin{proof}
We have \(S\simeq\mathbb P_{\prj{1}}(\shf{V})\), where \(\shf{V}=\sshf{\prj{1}}\oplus\sshf{\prj{1}}(e)\), and \(\sshf{S}(\xi)=\sshf{\mathbb P(\shf{V})}(1)\). By the projection formula,
\begin{equation}\label{eq:Fe-pushforward}
 \varphi_*\shf{G}
 \simeq\Sym^u\shf{V}\otimes\sshf{\prj{1}}(v)
 \simeq\bigoplus_{i=0}^u\sshf{\prj{1}}(v+ie).
\end{equation}
This bundle is globally generated, and the relative evaluation map \(\varphi^*\varphi_*\shf{G}\to\shf{G}\) is surjective since \(u\geq0\); hence \(\shf{G}\) is globally generated. As in the proof of \cref{lem:global-relative}\textup{(i)}, factoring the evaluation map of \(\shf{G}\) through \(\varphi^*\varphi_*\shf{G}\) gives an exact sequence
\[
0\longrightarrow\varphi^*M_{\varphi_*\shf{G}}
 \longrightarrow M_{\shf{G}}
 \longrightarrow M_{\shf{G}/\varphi}\longrightarrow0,
 \qquad M_{\shf{G}/\varphi}=\ker\bigl(\varphi^*\varphi_*\shf{G}\to\shf{G}\bigr).
\]
By \eqref{eq:Fe-pushforward}, \(\varphi^*M_{\varphi_*\shf{G}}\simeq\sshf{S}(-f)^{\oplus N_{\shf{G}}}\), since the kernel bundle of \(\sshf{\prj{1}}(d)\) is \(\sshf{\prj{1}}(-1)^{\oplus d}\).

It remains to compute the relative kernel bundle \(M_{\shf{G}/\varphi}\). If \(u=0\), the relative evaluation map is an isomorphism, and \(M_{\shf{G}/\varphi}=0\). Let \(u\geq1\), and consider the tautological sequence
\[
 0\longrightarrow\shf{K}\longrightarrow\varphi^*\shf{V}
 \longrightarrow\sshf{S}(\xi)\longrightarrow0,\qquad
 \shf{K}\simeq\varphi^*(\det\shf{V})\otimes\sshf{S}(-\xi)
 \simeq\sshf{S}(-\xi+ef).
\]
Multiplication in the symmetric algebra gives an exact sequence
\[
 0\longrightarrow
 \shf{K}\otimes\Sym^{u-1}(\varphi^*\shf{V})
 \longrightarrow\Sym^u(\varphi^*\shf{V})
 \longrightarrow\sshf{S}(u\xi)\longrightarrow0.
\]
Exactness is a local statement: where the tautological sequence splits, with \(\shf{K}\) generated by \(x\) and \(\varphi^*\shf{V}\) by \(x,y\), the symmetric algebra is \(\sshf{S}[x,y]\), the right-hand map sets \(x=0\), and its kernel in degree \(u\) is \(x\) times the degree-\((u-1)\) part. Since \(\varphi^*\varphi_*\shf{G}\simeq\Sym^u(\varphi^*\shf{V})\otimes\sshf{S}(vf)\) by \eqref{eq:Fe-pushforward}, twisting by \(\sshf{S}(vf)\) gives
\begin{align*}
M_{\shf{G}/\varphi} &\simeq \shf{K}\otimes\Sym^{u-1}(\varphi^*\shf{V})\otimes \sshf{S}(vf)
 \simeq\sshf{S}\bigl(-\xi+(v+e)f\bigr)
       \otimes\varphi^*\Sym^{u-1}\shf{V}\\
 &\simeq\bigoplus_{t=1}^u
       \sshf{S}\bigl(-\xi+(v+te)f\bigr).
\end{align*}
Substituting the two computations into the exact sequence above gives \eqref{eq:Fe-evaluation}.
\end{proof}

\begin{proposition}\label{prop:low-rank}
In the cyclic-cover setup, if \(\rho(S)\leq2\), then \(L\) satisfies
\(\NP_p\).
\end{proposition}

\begin{proof}
Suppose first that \(S=\prj{2}\), and let \(H\) be a line. Then \(A\sim aH\) with \(a=A\cdot H\geq p+3\), and \(\pi^*H\) is ample and basepoint-free. The adjoint theorem of Bangere--Lacini \cite[Theorem~1.2]{BangereLacini} gives Property \(\NP_p\) for \(K_X+a\pi^*H=L\).

It remains to consider \(S=\mathbb F_e\), with the notation of \cref{lem:hirzebruch-kernel}.

\smallskip\noindent\emph{Step 1: } Write
\[
 A\sim a\xi+cf,\qquad B\sim\alpha\xi+\delta f.
\]
Then \(a=A\cdot f\geq p+3\) and \(c=A\cdot C_0\geq p+3\), and \(\alpha=B\cdot f\geq0\) because \(rB\) is effective. Moreover,
\begin{equation}\label{eq:delta-dichotomy}
 \delta\geq0\qquad\text{or}\qquad\delta=-\frac er.
\end{equation}
Indeed, if \(\delta=B\cdot C_0<0\), then \cref{prop:fractional-negative} shows that \(C_0\) is a component of \(\Gamma\) with \(C_0^2=rB\cdot C_0\), that is, \(-e=r\delta\). Write \(D_i=a_i\xi+c_if\) for \(0\leq i<r\). Since \(K_S=-2\xi+(e-2)f\),
\begin{equation}\label{eq:ai-ci}
 a_i=a-2+(r-1-i)\alpha,\qquad
 c_i=c+e-2+(r-1-i)\delta.
\end{equation}
By \eqref{eq:delta-dichotomy} and \eqref{eq:ai-ci},
\begin{equation}\label{eq:ai-ci-bound}
 a_i,c_i\geq p+1\qquad(0\leq i<r).
\end{equation}
In particular, \(D=D_0\) has positive degree on \(f\) and on \(C_0\), and \(D^2>0\); hence \(D\) is ample, and so is \(L=\pi^*D\).

\smallskip\noindent\emph{Step 2:} Put
\[
 \shf{E}=\pi_*\sshf{X}(L)\simeq\bigoplus_{i=0}^{r-1}\sshf{S}(D_i).
\]
By \eqref{eq:ai-ci-bound} and \cref{lem:hirzebruch-kernel}, each \(\sshf{S}(D_i)\) is globally generated. Hence so is \(\shf{E}\), and \(M_{\shf{E}}\simeq\bigoplus_{i=0}^{r-1}M_{D_i}\). Since \(H^0(X,\sshf{X}(L))=H^0(S,\shf{E})\), the evaluation map of \(L\) factors as
\[
 H^0(S,\shf{E})\otimes\sshf{X}\longrightarrow\pi^*\shf{E}
 \longrightarrow\sshf{X}(L),
\]
where the second map is the natural map \(\pi^*\pi_*\sshf{X}(L)\to\sshf{X}(L)\). On the summand \(\pi^*\sshf{S}(D_0)=\sshf{X}(L)\) it is the identity. So it is surjective, and its kernel maps isomorphically, under the projection, onto \(\bigoplus_{i=1}^{r-1}\pi^*\sshf{S}(D_i)\). Thus \(L\) is globally generated, and the snake lemma gives an exact sequence
\begin{equation}\label{eq:cover-evaluation-Fe}
 0\longrightarrow\pi^*M_{\shf{E}}
 \longrightarrow M_L
 \longrightarrow
 \bigoplus_{i=1}^{r-1}\pi^*\sshf{S}(D_i)
 \longrightarrow0.
\end{equation}

\smallskip\noindent\emph{Step 3:} We claim that
\begin{equation}\label{eq:tensor-vanishing-Fe}
 H^1(X,M_L^{\otimes h}\otimes\sshf{X}(qL))=0
 \qquad(1\leq h\leq p+1,\ q\geq1).
\end{equation}
The filtration of \(M_L^{\otimes h}\) induced by \eqref{eq:cover-evaluation-Fe} has graded pieces \(\pi^*\bigl(M_{\shf{E}}^{\otimes u}\otimes\sshf{S}(D_{i_1}+\cdots+D_{i_v})\bigr)\), up to the order of the factors, where \(u+v=h\) and \(1\leq i_1,\dots,i_v\leq r-1\); for \(v=0\) the sum \(D_{i_1}+\cdots+D_{i_v}\) is zero. By the projection formula and \(\pi_*\sshf{X}\simeq\bigoplus_{j=0}^{r-1}\sshf{S}(-jB)\), it suffices to prove the vanishing of \(H^1\) for the bundles
\begin{equation}\label{eq:Fe-graded}
 M_{\shf{E}}^{\otimes u}\otimes
 \sshf{S}(D_{i_1}+\cdots+D_{i_v}+qD-jB)\qquad(0\leq j<r).
\end{equation}
Since \(qD-jB=(q-1)D+D_j\), the line bundle in \eqref{eq:Fe-graded} is \(\sshf{S}(g\xi+kf)\) with \(g,k\geq p+1\), by \eqref{eq:ai-ci-bound}. By \eqref{eq:Fe-evaluation}, each summand \(M_{D_i}\) of \(M_{\shf{E}}\) has a filtration with graded pieces \(\sshf{S}(-f)\) and \(\sshf{S}(-\xi+(c_i+te)f)\), \(t\geq1\), where \(c_i+te\geq p+1\). Hence \eqref{eq:Fe-graded} has a filtration whose graded pieces are line bundles
\begin{equation*}
 \sshf{S}\bigl((g-s)\xi+(k+\gamma-(u-s))f\bigr),
 \qquad 0\leq s\leq u,\quad \gamma\geq s(p+1).
\end{equation*}
Since \(s\leq u\leq h\leq p+1\), both coefficients are nonnegative. For \(x,y\geq0\), the sheaf \(\varphi_*\sshf{S}(x\xi+yf)\simeq\Sym^x\shf{V}\otimes\sshf{\prj{1}}(y)\) is a direct sum of line bundles of nonnegative degree, and \(R^1\varphi_*\sshf{S}(x\xi+yf)=0\); hence \(H^1(S,\sshf{S}(x\xi+yf))=0\). This proves \eqref{eq:tensor-vanishing-Fe}.

In characteristic zero, \(\bigwedge^hM_L\) is a direct summand of \(M_L^{\otimes h}\), so \eqref{eq:tensor-vanishing-Fe} gives \[H^1\bigl(X,\bigwedge^hM_L\otimes\sshf{X}(qL)\bigr)=0\] for \(1\leq h\leq p+1\) and \(q\geq1\). By the kernel-bundle description preceding
\cref{lem:kernel-criterion}, these vanishings imply
\(K_{a,q}(X,L)=0\) for \(0\leq a\leq p\) and \(q\geq2\). Since \(L\) is very ample by \cref{lem:adjoint-separation}\textup{(ii)}, it satisfies \(\NP_p\).
\end{proof}

\begin{proof}[Proof of \cref{thm:main}]
If \(\rho(S)\leq2\), the result follows from \cref{prop:low-rank}, so
assume that \(\rho(S)\geq3\). If \(q(X)=0\), apply
\cref{prop:regular-case}. We may therefore assume that \(q(X)>0\).
By \cref{lem:regularity}, this forces \(P_B^2=0\).

If \(\dim|-K_S|\geq1\), \cref{cor:moving-anti} proves the result.
Suppose that \(h^0(S,\sshf{S}(-K_S))=1\). Then \(K_S^2\leq0\) by
\cref{lem:anti-subcurves}\textup{(iii)}. If \(K_S^2=0\), \cref{thm:nonnegative} applies. Otherwise \(K_S^2<0\), so \eqref{eq:fixed-boundary} holds, and \cref{thm:fixed-uniform} applies. This completes the proof.
\end{proof}

{\small\noindent\textsc{School of Mathematics, Sun Yat-sen University\\
W. 135 Xingang Rd., Guangzhou, Guangdong 510275, P.R.~China}\\
\emph{E-mail address}: \texttt{songlei3@mail.sysu.edu.cn}\par}

\end{document}